\documentclass[11pt]{article}
\usepackage[margin=1in]{geometry}
\usepackage{amsmath,amssymb,amsthm,mathtools}
\usepackage{enumitem}
\usepackage{booktabs,array}
\usepackage{placeins}
\usepackage[authoryear,round]{natbib}
\usepackage{xcolor}
\usepackage{multirow}
\usepackage[colorlinks=true,citecolor=blue!70!black,linkcolor=red!60!black,urlcolor=blue!70!black]{hyperref}
\usepackage{authblk}

\newtheorem{theorem}{Theorem}
\newtheorem{proposition}{Proposition}
\newtheorem{lemma}{Lemma}

\theoremstyle{definition}

\theoremstyle{remark}

\def\bbR{\mathbb{R}}
\newcommand{\E}{\mathbb{E}}
\renewcommand{\P}{\mathbb{P}}
\newcommand{\Var}{\operatorname{Var}}
\newcommand{\Cov}{\operatorname{Cov}}
\newcommand{\FDR}{\mathrm{FDR}}
\newcommand{\FDP}{\mathrm{FDP}}
\newcommand{\Ntil}{\widetilde N}
\newcommand{\Stil}{\widetilde S}
\newcommand{\logp}{\log_{+}}

\title{A Dimension-Free Bound on the FDR of the Benjamini–Hochberg Procedure under Slowly Decaying Gaussian Correlations}

\date{\today}
\author[1]{Zuofeng Shang}
\author[1]{Wenge Guo}
\affil[1]{Department of Mathematical Sciences, New Jersey Institute of Technology}

\begin{document}
\maketitle

\begin{abstract}
The Benjamini–Hochberg (BH) procedure is proven to control the false discovery rate under independence and under positive dependence. Under arbitrary dependence, however, its guarantee deteriorates by a factor that grows logarithmically with the number of hypotheses. For Gaussian test statistics, it is unknown whether dependence can make the false discovery rate grow with the number of tests. We study this question when the correlation between two test statistics decays with the distance between their indices. We seek guarantees that hold uniformly over a wide class of correlation structures and over all configurations of the signals. We show that if the correlations decay at least as fast as a power of the inverse logarithm of the distance with exponent larger than one, then the false discovery rate of the original BH procedure is bounded by the nominal level times a factor that does not depend on the number of hypotheses. This factor grows only polylogarithmically as the level decreases. The result holds for both one- and two-sided tests, and the procedure uses no knowledge of the correlations. The class allows correlations of both signs. The required decay can be arbitrarily close to the inverse logarithm. The key technical tool is a truncated Hermite expansion of the exceedance indicators. It shows that their main fluctuations are common to all thresholds that BH may use. These fluctuations can be well controlled, while the remainder is too small to accumulate which leads to a bound that does not inflate with the dimension.
\\

\end{abstract}


\section{Introduction}
Modern scientific studies routinely test thousands or millions of hypotheses at once, and the test statistics involved are rarely independent. For instance, expression levels of co-regulated genes move together, neighbouring loci along a genome and adjacent voxels in a brain image are spatially correlated, and statistics computed along a time series inherit its serial dependence. The BH procedure of \citet{BH1995} has become the standard tool for controlling the false discovery rate (FDR) in such studies. It is natural to seek guarantees when independence fails.

The original analysis of \citet{BH1995} shows that, for independent test statistics, the FDR of the procedure equals the nominal level multiplied by the proportion of true null hypotheses. \citet{BY2001} and \citet{Sarkar2002} extended this guarantee to positive regression dependence, a condition that for one-sided Gaussian tests amounts to nonnegative correlations but is much harder to verify for two-sided tests \citep{KarlinRinott1981}. Without any assumption on the dependence, \citet{BY2001} showed that the FDR can exceed the nominal level only by a factor of the order of the logarithm of the number of hypotheses. This factor cannot be improved over all joint distributions, but the distributions that attain it are far from Gaussian, and in practice the procedure is rarely observed to behave so badly. A large body of work has therefore studied weak dependence, showing that under mixing or decay conditions the false discovery proportion behaves asymptotically as under independence \citep{Farcomeni2007,Wu2008,ClarkeHall2009,DelattreRoquain2016,DattaDey2026}, and that dependence mainly affects its variability \citep{Owen2005,Efron2007,SchwartzmanLin2011}. These results describe the limit as the number of hypotheses grows for a fixed dependence structure. 
They do not provide bounds on the FDR that hold for every number of hypotheses and uniformly over the possible correlations and signals, which is precisely a guarantee under unknown dependence.

Several recent works obtain bounds on the FDR of BH that do not grow with $m$ under conditions on the null $p$-values only. \citet{Su2018} considered positive regression dependence among the null $p$-values, and \citet{ChiRamdasWang2025} considered weak negative dependence. For Gaussian test statistics, this question has been taken up from several directions. First, Gaussian dependence can inflate the FDR. \citet{Dobriban2026} showed that BH can fail to control the FDR at the nominal level when two-sided Gaussian tests are correlated. \citet{Lei2026} showed that the inflation can exceed any fixed multiple of the nominal level, and conjectured that the worst-case FDR does not grow with $m$. Second, bounds free of $m$ have been obtained for special structures. \citet{Lei2026} proved such bounds for common-factor models. \citet{Zhang2026} established exact Simes inequalities for two-sided Gaussian tests, which control the FDR under the global null for several classes of covariance matrices. \citet{LyuWang2026} studied one-sided tests in fixed dimensions, and obtained bounds free of $m$ under the global null with nonpositive correlations and in rank-one models. Other recent work derived finite-sample bounds on the FDR of BH for a given covariance matrix \citep{GhoshSarkar2026}, proposed shifted procedures for two-sided tests \citep{Sarkar2023,GhoshSarkar2025}, and restored exact control when the dependence is known \citep{FithianLei2022,NguyenFithian2025}.
To our knowledge, it has remained open whether the logarithmic factor of \citet{BY2001} is genuinely needed for the original BH procedure uniformly over a wide class of Gaussian correlation structures that allows correlations to change sign and to decay slowly, when the signal strengths and their locations are arbitrary.

In this short note, we answer this question when the correlation between two test statistics decays with the distance between their indices. This includes stationary dependence, in which the correlation depends only on this distance. We show that if the absolute correlations decay at least as fast as a power of the inverse logarithm of this distance with exponent larger than one, then the FDR of the BH procedure is bounded by the nominal level times a factor that does not depend on the number of hypotheses. No condition is imposed on the signs of the correlations, which may change arbitrarily. The bound holds uniformly over all such correlation matrices and all configurations of the signals, for one- and two-sided tests, and the required decay is arbitrarily close to the inverse logarithm. A main challenge is that the threshold of BH is random, so the exceedances must be controlled at all possible thresholds. A truncated Hermite expansion shows that their main fluctuations are common to all thresholds and can be easily controlled, while the remainder of the Hermite expansion is too small to accumulate. This yields an FDR bound that does not inflate with the number of hypotheses.

The rest of the paper is organized as follows. Section~\ref{sec:setting} introduces the model, defines the worst-case FDR, states the main result, discusses the dependence scope of our result, and sketches its proof. Section~\ref{sec:numerics} compares the FDR of BH with the BY bound and our bound numerically. Section~\ref{sec:appendix} contains the proofs of main and preliminary results.

\section{Multiple Testing and FDR}\label{sec:setting}
Suppose we observe a Gaussian vector $T=(T_1,\ldots,T_m)^\top\sim N(\mu,\Sigma)$ with unknown mean vector $\mu=(\mu_1,\ldots,\mu_m)^\top$ and unknown covariance matrix $\Sigma$ whose diagonal entries equal one. We test the null hypotheses
\[
H_i^0:\mu_i=0,\qquad i=1,\ldots,m,
\]
against either one-sided alternatives $H_i^{1,a}:\mu_i>0$ or two-sided alternatives $H_i^{2,a}:\mu_i\ne0$. Let $\mathcal H_0=\{i:\mu_i=0\}$ be the set of true null hypotheses, $m_0=|\mathcal H_0|$ be the number of true null hypotheses, and $\pi_0=m_0/m$ be the proportion of true null hypotheses. 
With $\Phi$ the standard normal cumulative distribution function, the one-sided and two-sided $p$-values are
\begin{eqnarray*}
(\textrm{one-sided}):&& P_i^1=1-\Phi(T_i),\\
(\textrm{two-sided}):&& P_i^2=2\big(1-\Phi(|T_i|)\big),
\end{eqnarray*}
both uniformly distributed on $[0,1]$ under $H_i^0$. Throughout, $P_i$ denotes either $P_i^1$ for all $i$ (one-sided case) or $P_i^2$ for all $i$ (two-sided case). Let $P_{(1)}\le\cdots\le P_{(m)}$ be the ordered $p$-values. At level $\alpha\in(0,1)$, BH rejects $H_i^0$ if and only if $P_i\le\alpha R/m$, where
\[
R=\max\{1\le i\le m:P_{(i)}\le\alpha i/m\}
\]
is the number of rejections, with $R=0$ (no rejection) if the set is empty. Let $V$ be the number of rejected true null hypotheses. The FDR is
\begin{equation}\label{def:FDR}
\FDR(\alpha;m,\mu,\Sigma)=\E\Big[\frac{V}{R\vee1}\Big],
\end{equation}
where $a\vee b=\max\{a,b\}$.
The notation suppresses the dependence on whether one-sided or two-sided $p$-values are used. All our results hold in both cases.

\subsection{Worst-Case FDR under Decaying Correlations}
Let $\mathbb N_0=\{0,1,2,\ldots\}$. For $\epsilon>0$ and $C>0$, define the envelope $\rho_{\epsilon,C}\colon\mathbb N_0\to[0,1]$ by $\rho_{\epsilon,C}(0)=1$ and
\[
\rho_{\epsilon,C}(k)=\min\big\{1,\ C\log^{-(1+\epsilon)}(k+e)\big\},\qquad k\ge1.
\]
Consider the class of correlation matrices
\[
\mathcal C_m(\epsilon,C)=\big\{\Sigma\in\bbR^{m\times m}:\Sigma=\Sigma^\top,\ \Sigma\succeq0,\ \Sigma_{ii}=1,\ |\Sigma_{ij}|\le\rho_{\epsilon,C}(|i-j|)\text{ for all }i\ne j\big\}.
\]
Thus the correlation between two test statistics is bounded by the envelope evaluated at the distance between their indices. 
An important special case is stationary dependence, where $\Sigma_{ij}=\rho(|i-j|)$ for a correlation sequence $\rho$. Then the condition reads $|\rho(k)|\le\rho_{\epsilon,C}(k)$ for $1\le k\le m-1$.
Since $\rho_{\epsilon,C}(k)\log k\to0$, every such $\rho$ satisfies Berman's condition $\rho(k)\log k\to0$ \citep{Berman1964}. 
As $\epsilon$ approaches zero, the envelope approaches the critical rate $\log^{-1}k$ of \citet{MittalYlvisaker1975}. To study the FDR of BH under such slowly decaying correlations, we consider the worst-case FDR over this class:
\[
U(\alpha;\epsilon,C)=\sup_{m\ge1}\ \sup_{\Sigma\in\mathcal C_m(\epsilon,C)}\ \sup_{\mu\in\bbR^m}\ \FDR(\alpha;m,\mu,\Sigma).
\]
For arbitrary $\mu$ and $\Sigma$, \citet{BY2001} showed that
\begin{equation}\label{BH:original:bound}
\FDR(\alpha;m,\mu,\Sigma)\le\alpha\pi_0H_m,
\end{equation}
where $H_m\triangleq\sum_{k=1}^mk^{-1}\sim\log m$. This bound grows with $m$ and is sharp under arbitrary dependence among the $T_i$, but the distributions known to attain it are non-Gaussian. It is unclear whether a bound free of $m$ holds when $\Sigma$ is restricted to $\mathcal C_m(\epsilon,C)$, that is, whether $U(\alpha;\epsilon,C)<\infty$.

\subsection{Main Result}\label{sec:main}
In this section, we establish an upper bound for the worst-case FDR over the mean vector $\mu$ and the class of correlation matrices $\mathcal C_m(\epsilon,C)$, that does not depend on the number of hypotheses $m$.
Let $\mathcal S=\{1,2,3,\dots\}$ for one-sided $p$-values and $\mathcal S=\{2,4,6,\dots\}$ for two-sided $p$-values, and fix an element $k^*\in\mathcal S$. Let $\phi$ be the standard normal density and $He_k$ the probabilists' Hermite polynomials, defined by 
\[
\textrm{$He_0(x)\equiv 1$,\,\,$\phi^{(k)}(x)=(-1)^kHe_k(x)\phi(x)$ for $k=1,2,\ldots$.}
\] 
For $w\in(0,1/2]$, let 
\[
z_w=\left\{\begin{array}{ll}\Phi^{-1}(1-w),&\textrm{(one-sided)},\\[2pt] 
\Phi^{-1}(1-w/2),&\textrm{(two-sided)}.\end{array}\right.
\]
Here $z_w$ is the upper $w$-quantile of $\Phi$ in the one-sided case and the upper $w/2$-quantile of $\Phi$ in the two-sided case.
For $z\ge0$ and $r\in[-1,1]$, let 
\[
\textrm{$r_z=\min\{1/2,z^{-2}\}$ and $R(r,z)=\textrm{Cov}(\mathbf 1\{X>z\}, \mathbf 1\{Y>z\})$,}
\]
where $(X,Y)$ is the standard bivariate normal pair with correlation coefficient $r$.
Here $R(r,z)$ is the covariance of the two exceedance indicators $\mathbf 1\{X>z\}$ and $\mathbf 1\{Y>z\})$. Let
\[
\tau(r,z)=\begin{cases}2,&|r|>1/2,\\ 2\min\left\{1,\frac{R(|r|,z)}{1-\Phi(z)}\right\},&|r|\le1/2.\end{cases}
\]
Let $J=\lceil\log_2m\rceil$, $w_j=2^{-j}\alpha$ for $j=0,\ldots,J-1$, and $z_j=z_{w_j}$.
For $k\ge1$, let 
\[
G_{\epsilon,C,k}(m)=\sum_{l=1}^{m-1}\rho_{\epsilon,C}(l)^k,\,\,\,\,
b_k=\sup_{0\le z\le z_{J-1}}\frac{\phi(z)|He_{k-1}(z)|}{k!\,(1-\Phi(z))}.
\]
Here $G_{\epsilon,C,k}(m)$ is the sum of the $k$-th powers of the envelope $\rho_{\epsilon,C}(l)$ over the lags $l=1,\ldots,m-1$. It can be shown that the exceedance indicator $\mathbf 1\{X>z\}$ has $k$-th Hermite coefficient $\phi(z)He_{k-1}(z)/k!$; see Lemma \ref{lem:hermite}. Thus $b_k$ bounds this coefficient, relative to the mean $1-\Phi(z)$ of the indicator, uniformly over the thresholds $z_0,\ldots,z_{J-1}$.
For $z\ge0$, let
\begin{eqnarray*}
\mathsf A(z)&=&\sum_{k\in\mathcal S,\ k\ge k^*}\frac{|He_{k-1}(z)|^2}{k!}\,r_z^{\,k-k^*},\,\,\,\,
e_{\epsilon,C,m}(z)=2\Big(\frac{\phi(z)}{1-\Phi(z)}\Big)^2\mathsf A(z)\,\frac{G_{\epsilon,C,k^*}(m)}m.
\end{eqnarray*}
The function $\mathsf A(z)$ collects the squared Hermite coefficients of order $k\ge k^*$ of the exceedance indicator at the threshold $z$, normalized by $\phi(z)^2$ and powers of $r_z$. The function $e_{\epsilon,C,m}(z)$ combines it with the average $G_{\epsilon,C,k^*}(m)/m$ of the weak correlations over the lags, and bounds their contribution to the relative variance of the exceedance count at $z$. 
Define
\begin{equation}\label{eq:K}
K_{\epsilon,C}=1+\sup_{z\ge0}\sum_{k\ge1:\,\rho_{\epsilon,C}(k)>r_z}\tau\big(\rho_{\epsilon,C}(k),z\big),
\end{equation}
and
\begin{equation}\label{eq:Lambda:epsilon:C:m}
\Lambda_{\epsilon,C,m}(\alpha)=\left\{\begin{array}{ll}\Big(\sum_{k\in\mathcal S,\ k<k^*}b_k\Big[k!\Big(\frac1m+\frac{2G_{\epsilon,C,k}(m)}m\Big)\Big]^{1/2}\Big)^2, & k^*>\min\mathcal S,\\
0, & k^*=\min\mathcal S.\end{array}\right.
\end{equation}
The constant $K_{\epsilon,C}$ collects the strongly correlated lags, and the constant $\Lambda_{\epsilon,C,m}(\alpha)$ involves the Hermite terms of order $k<k^*$. The former accounts for the few lags with large correlations. The latter accounts for a random fluctuation that is common to all thresholds $z_0,\ldots,z_{J-1}$, up to deterministic coefficients. The coefficients are bounded uniformly over these thresholds through $b_k$, so this fluctuation is controlled once for all of them and does not accumulate over thresholds.

Throughout, we fix $\alpha\in(0,1/2]$. For $K\ge1$, $E\ge0$ and $\pi_0\in(0,1]$, let
\[
g(t;K,E,\pi_0)=\frac{6\pi_0\sqrt K}{t-2\pi_0}+\frac{4\pi_0E}{(t-2\pi_0)^2},\,\,\,t>2\pi_0,
\]
and
\begin{equation}\label{eq:def:Psi}
\Psi_{\pi_0}(K,E,\alpha)=2\pi_0+\int_{2\pi_0}^{1/\alpha}\min\{1,g(t;K,E,\pi_0)\}\,dt.
\end{equation}
Since $\pi_0\le1$ and $\alpha\le1/2$, we have $2\pi_0\le1/\alpha$, so the above integral limits are valid.
Letting $\pi_0\to0$ in (\ref{eq:def:Psi}), we get $\Psi_0(K,E,\alpha)=0$. 
For $m\ge2$, let 
\[
E_{\epsilon,C,m}(\alpha)=\sum_{j=0}^{J-1}e_{\epsilon,C,m}(z_j).
\]
The quantity $E_{\epsilon,C,m}(\alpha)$ adds $e_{\epsilon,C,m}$ over the thresholds $z_0,\ldots,z_{J-1}$. Since the threshold of BH is random, the exceedance counts must be controlled at all of these thresholds simultaneously. Thus $E_{\epsilon,C,m}(\alpha)$ measures the total effect of the weakly correlated lags on the FDR. When $\epsilon k^*\ge1$, it stays bounded in $m$. Define
\begin{equation}\label{eq:Fm}
F_{\epsilon,C,m}(\alpha;\pi_0)=\inf_{\theta\in(0,1)}\Psi_{\pi_0}\Big(\frac{K_{\epsilon,C}}{(1-\theta)^2},\ \frac{\Lambda_{\epsilon,C,m}(\alpha)}{\theta^2}+\frac{E_{\epsilon,C,m}(\alpha)}{(1-\theta)^2},\ \alpha\Big).
\end{equation}
When $m=1$, let $F_{\epsilon,C,1}(\alpha;\pi_0)=\pi_0$. If $k^*=\min\mathcal S$, then $\Lambda_{\epsilon,C,m}(\alpha)=0$, and hence, (\ref{eq:Fm}) becomes
\begin{equation}\label{eq:computing:F}
F_{\epsilon,C,m}(\alpha;\pi_0)=\Psi_{\pi_0}(K_{\epsilon,C},E_{\epsilon,C,m}(\alpha),\alpha).
\end{equation}
The following theorem establishes an upper bound for the FDR of HG procedure.
\begin{theorem}\label{thm:main}
Let $\alpha\in(0,1/2]$, $\pi_0\in[0,1]$, $\epsilon>0$, $C>0$ and $k^*\in\mathcal S$.
\begin{enumerate}[label=(\roman*).,ref=(\roman*)]
\item\label{thm:main:i} For both one-sided and two-sided $p$-values, for any $m\ge1$, $\mu\in\bbR^m$ and $\Sigma\in\mathcal{C}_m(\epsilon,C)$,
\begin{equation}\label{our:FDR:bound}
\FDR(\alpha;m,\mu,\Sigma)\le\alpha\min\Big\{\pi_0H_m,\ F_{\epsilon,C,m}(\alpha;\pi_0)\Big\}.
\end{equation}

\item\label{thm:main:ii} Moreover, if $\epsilon k^*\ge1$, then 
\begin{equation}\label{eq:supF}
\sup_{m\ge2}F_{\epsilon,C,m}(\alpha;\pi_0)\le A(\epsilon,C,\pi_0)+B(\epsilon,C,\pi_0)\log^{\max\{1,k^*/2\}}(1/\alpha),
\end{equation}
where, for $\pi_0>0$,
\begin{align*}
A(\epsilon,C,\pi_0)={}&\pi_0\big(2+12\sqrt{\bar K_{C,k^*}}\big)+12\pi_0\sqrt{\bar K_{C,k^*}}\Big(1+\logp\frac1{12\pi_0\sqrt{\bar K_{C,k^*}}}\Big)\\
&+4\max\{2,2^{k^*/2}\}\sqrt{\pi_0\big(\bar\Lambda_{C,k^*}+\bar E_{C,k^*}\big)},\\
B(\epsilon,C,\pi_0)={}&12\pi_0\sqrt{\bar K_{C,k^*}}+4\max\{2,2^{k^*/2}\}\sqrt{\pi_0\big(\bar\Lambda_{C,k^*}+\bar E_{C,k^*}\big)},
\end{align*}
and $A(\epsilon,C,0)=B(\epsilon,C,0)=0$. Here $\logp x=\max\{\log x,0\}$. With $p=1+1/k^*$ and $N=\#\{k\in\mathcal S:k<k^*\}$, the constants 
$\bar K_{C,k^*}$, $\bar\Lambda_{C,k^*}$ and $\bar E_{C,k^*}$ are given as follows
\begin{align}
\bar K_{C,k^*}&=1+2\big(e^{(2C)^{1/p}}-e\big)_++\frac{2}{\sqrt{6\pi}}\sup_{z>\sqrt2}(1+z)e^{-z^2/6}\Big(e^{(Cz^2)^{1/p}}-e\Big)_+,\label{eq:Kbar}\\
\bar E_{C,k^*}&=a_*C^{k^*}\Big(\frac4{\log2}\Big)^{k^*+1}\Big[\Big(\frac{2k^*+2}e\Big)^{k^*+1}+2^{k^*+1}\Big],\label{eq:Ebar}\\
\bar\Lambda_{C,k^*}&=\frac N2\sum_{k\in\mathcal S,\ k<k^*}\frac1{k!}\Big(\frac{8k+8}{\log2}\Big)^k\Big[\Big(\frac ke\Big)^k+2C^k\Big(\Big(\frac{2k}e\Big)^k+\frac{2^{kp}}{\log2}\Big)\Big],\label{eq:Lbar}
\end{align}
where $a_*=2/\sqrt3$ in the two-sided case with $k^*=2$, and $a_*=2(e-1)/\sqrt3$ otherwise. If $k^*=\min\mathcal S$, then $N=0$ and hence one can set $\bar\Lambda_{C,k^*}=0$.
\end{enumerate}
\end{theorem}
The proof of Theorem \ref{thm:main} is given in Section~\ref{sec:appendix}.
Theorem~\ref{thm:main}\ref{thm:main:i} gives an upper bound on the FDR of the BH procedure under decaying Gaussian correlations. It holds for every $m\ge1$, $\mu\in\bbR^m$ and $\Sigma\in\mathcal{C}_m(\epsilon,C)$. The bound is the smaller of the BY bound $\alpha\pi_0H_m$ and $\alpha F_{\epsilon,C,m}(\alpha;\pi_0)$, so it never exceeds the BY bound. Given $m$, $\epsilon$, $C$, $\pi_0$ and $\alpha$, exact calculation of the quantity $F_{\epsilon,C,m}(\alpha;\pi_0)$ is hard, but numerical computation is feasible.
Table~\ref{tab:F2} reports numerical values of $F_{\epsilon,C,m}(\alpha;\pi_0)$ for $m=10^4$ and $m=10^8$ under several values of $\epsilon$, $C$, $k^*$, $\pi_0$ and $\alpha$. For the values of $\epsilon$ and $C$ being considered, the values appear to change little. The values decrease quickly as $\pi_0$ decreases. The values appear to change little with $m$, which illustrates that the FDR bound is free of $m$.

\begin{table}[htp]
\centering\small
\caption{\it Values of $F_{\epsilon,C,m}(\alpha;\pi_0)$ under various situations with $m=10^4$ and $m=10^8$.}\label{tab:F2}
\begin{tabular}{llccccccc}
\toprule
& & & \multicolumn{3}{c}{$m=10^4$} & \multicolumn{3}{c}{$m=10^8$}\\
\cmidrule(lr){4-6}\cmidrule(lr){7-9}
$p$-values & $(\epsilon,C,k^*)$ & $\pi_0$ & $\alpha=0.01$ & $\alpha=0.05$ & $\alpha=0.1$ & $\alpha=0.01$ & $\alpha=0.05$ & $\alpha=0.1$\\
\midrule
\multirow{4}{*}{One-sided} & \multirow{4}{*}{$(0.5,0.01,2)$} & 1 & 27.79 & 15.92 & 9.98 & 27.14 & 15.66 & 9.96\\
 &  & 0.5 & 16.73 & 10.85 & 8.10 & 16.25 & 10.60 & 7.96\\
 &  & 0.2 & 8.33 & 5.94 & 4.86 & 8.02 & 5.74 & 4.72\\
 &  & 0.1 & 4.85 & 3.63 & 3.08 & 4.64 & 3.48 & 2.97\\
\midrule
\multirow{4}{*}{One-sided} & \multirow{4}{*}{$(1,0.01,1)$} & 1 & 24.83 & 14.63 & 9.74 & 24.81 & 14.62 & 9.74\\
 &  & 0.5 & 14.56 & 9.59 & 7.33 & 14.54 & 9.57 & 7.32\\
 &  & 0.2 & 6.98 & 5.01 & 4.14 & 6.95 & 4.99 & 4.13\\
 &  & 0.1 & 3.94 & 2.96 & 2.53 & 3.92 & 2.94 & 2.51\\
\midrule
\multirow{4}{*}{Two-sided} & \multirow{4}{*}{$(1,0.1,2)$} & 1 & 24.77 & 14.60 & 9.73 & 24.76 & 14.59 & 9.73\\
 &  & 0.5 & 14.50 & 9.54 & 7.30 & 14.49 & 9.54 & 7.30\\
 &  & 0.2 & 6.91 & 4.96 & 4.10 & 6.90 & 4.95 & 4.10\\
 &  & 0.1 & 3.88 & 2.91 & 2.48 & 3.87 & 2.90 & 2.48\\
\bottomrule
\end{tabular}
\end{table}
Theorem~\ref{thm:main}\ref{thm:main:ii} describes how $F_{\epsilon,C,m}(\alpha;\pi_0)$ grows with $\log(1/\alpha)$, uniformly in $m$. This is useful for understanding its behavior as $\alpha\to0$. Specifically, together with \eqref{our:FDR:bound}, it shows that the FDR of BH is at most of order $\alpha\log^{\max\{1,k^*/2\}}(1/\alpha)$, uniformly over $m$, $\mu$ and $\Sigma\in\mathcal C_m(\epsilon,C)$. In particular, it tends to zero as $\alpha\to0$. When $\epsilon\ge1/2$, one may take $k^*\le2$, and the order becomes $\alpha\log(1/\alpha)$.
It should be emphasized that the constants $A(\epsilon,C,\pi_0)$ and $B(\epsilon,C,\pi_0)$ are only for demonstration purposes but are conservative. In practice, we suggest to use (\ref{eq:computing:F}) to compute the value of $F_{\epsilon,C,m}(\alpha;\pi_0)$.

Classical FDR bounds either assume positive regression dependence or grow like $\log m$ \citep{BY2001,Sarkar2002}, and results under weak dependence are asymptotic in $m$ \citep{Farcomeni2007,Wu2008,ClarkeHall2009,DelattreRoquain2016}. Table~\ref{tab:dimension-free-fdr} compares Theorem~\ref{thm:main} with recent bounds that do not grow with $m$, with rates stated as $\alpha\to0$.

\begin{table}[htp]
\centering\small
\caption{\it Upper bounds on the FDR of BH that do not depend on $m$.}\label{tab:dimension-free-fdr}
\setlength{\tabcolsep}{4pt}
\renewcommand{\arraystretch}{1.15}
\begin{tabular}{@{}p{0.2\linewidth}p{0.5\linewidth}p{0.25\linewidth}@{}}
\toprule
Result & Setting & FDR upper bound\\
\midrule
\citet{Su2018} & Positive regression dependence among the true-null $p$-values. & $O(\alpha\log(1/\alpha))$\\
\citet{ChiRamdasWang2025} & Weak negative dependence among the true-null $p$-values. & $O(\alpha\log(1/\alpha))$\\
\citet{Zhang2026} & Two-sided Gaussian tests with bounded covariance degree among the true-null statistics. & $O(\alpha\log(1/\alpha))$\\
\citet{Lei2026} & Two-sided Gaussian tests in a class of common-factor models. & $O(\alpha\sqrt{\log(1/\alpha)})$\\
\citet{LyuWang2026} & One-sided Gaussian tests with rank-one dependence and monotone non-null $p$-value maps. & $2W(1)\alpha$ for $\alpha<1/2$\\
Theorem~\ref{thm:main} & One- and two-sided Gaussian tests with $\Sigma\in\mathcal{C}_m(\epsilon,C)$ for some $\epsilon,C>0$, where $k^*$ is the smallest element of $\mathcal S$ with $\epsilon k^*\ge1$. & $O\big(\alpha\log^{\max\{1,k^*/2\}}(1/\alpha)\big)$\\
\bottomrule
\end{tabular}
\end{table}

\subsection{Proof Sketch}\label{sec:sketch}
Let $\xi=T-\mu$ and let $\Ntil(w)$ be the number of true nulls whose $p$-values are at most $w$. Every rejected true null is counted by $\Ntil$ at the BH threshold, so the false discovery proportion is bounded by the supremum of $\Ntil(u)/(mu)$ over the range of thresholds used by BH. This supremum depends on $\mu$ only through the set of true nulls (Lemma~\ref{lem:reduction}), and the FDR is bounded by an integral of its tail probabilities. The main goal of the proof is therefore a tail bound for this normalized supremum that does not depend on $m$. Integrating it then gives a bound on the FDR of the form $\alpha$ times a constant free of $m$.

To obtain the tail bound, we first restrict the thresholds to dyadic levels, which costs a factor $2$ (Lemma~\ref{lem:dyadic}). At each level we expand the exceedance indicators generated by $\xi$ in Hermite polynomials and split the centered count into low-order terms and a remainder (Lemma \ref{lem:hermite}). The low-order terms are, up to deterministic coefficients, the same random variables at every level, so they are controlled once rather than once per level, which avoids a union bound over the logarithmically many levels. The remainder is controlled by Chebyshev's inequality at coarse levels and by a first-moment bound at fine levels. Summing the better of the two over the levels costs only a constant factor (Lemma~\ref{lem:geom}). The variance of the remainder (Proposition~\ref{prop:var}) has a local part, coming from the finitely many strongly and moderately correlated lags, and a long-range part, coming from the many weakly correlated lags. Both parts depend on the correlations only through their absolute values and increase with them, so it suffices to evaluate them at the envelope. Summed over the levels, the long-range part stays bounded in $m$ when $\epsilon k^*\ge1$ (Section~\ref{sec:proofmain}). It grows polylogarithmically in $1/\alpha$, which gives the rate in Theorem~\ref{thm:main}\ref{thm:main:ii}. When $\epsilon\ge1$ for one-sided tests or $\epsilon\ge1/2$ for two-sided tests, one may take $k^*=\min\mathcal S$. Then no low-order terms are needed and the argument reduces to a plain second-moment bound.

\subsection{Scope of the Dependence Condition}\label{subsec:dependence-scope}
In this section, we discuss several Gaussian dependence models that satisfy the dependence conditions of Theorem~\ref{thm:main}, as well as several models that do not. For each model, we examine its correlations and then check whether its correlation matrix lies in $\mathcal C_m(\epsilon,C)$. For a stationary model with correlation sequence $\rho$, this holds whenever $|\rho(k)|\le\rho_{\epsilon,C}(k)$ for all $k\ge1$. Since $\Sigma$ is the correlation matrix of $T-\mu$, we describe each model through $\xi_i=T_i-\mu_i$. Let $Z_i,i\in\mathbb{Z}$ be independent standard normal variables.
The following Table~\ref{tab:scope} summarizes the models.
\begin{table}[htp]
\centering\small
\caption{\it Gaussian dependence models and the class $\mathcal C_m(\epsilon,C)$.}\label{tab:scope}
\setlength{\tabcolsep}{4pt}
\renewcommand{\arraystretch}{1.15}
\begin{tabular}{@{}>{\raggedright\arraybackslash}p{0.24\linewidth}>{\raggedright\arraybackslash}p{0.27\linewidth}>{\raggedright\arraybackslash}p{0.08\linewidth}>{\raggedright\arraybackslash}p{0.36\linewidth}@{}}
\toprule
Model & Correlations & Covered & Choice of $(\epsilon,C)$\\
\midrule
AR(1), $|a|<1$ & $\rho(k)=a^k$ & Yes & $\epsilon>0$, $C=\sup_{k\ge1}|a|^k\log^{1+\epsilon}(k+e)$\\
MA($q$) & $\rho(k)=0$ for $k>q$ & Yes & $\epsilon>0$, $C=\log^{1+\epsilon}(q+e)$\\
Fractional Gaussian noise, $H\in(1/2,1)$ & $\rho(k)\sim H(2H-1)k^{2H-2}$ & Yes & $\epsilon>0$, a suitable $C$\\
Independent blocks of size at most $B$ & $\Sigma_{ij}=0$ for $|i-j|\ge B$ & Yes & $\epsilon>0$, $C=\log^{1+\epsilon}(B-1+e)$\\
\midrule
Equicorrelation, fixed $r>0$ & $\Sigma_{ij}=r$ for all $i\ne j$ & No & None, since $\Sigma_{1m}=r$ while $\rho_{\epsilon,C}(m-1)\to0$\\
One-factor, $|\lambda_i|\ge\lambda>0$ & $\Sigma_{ij}=\lambda_i\lambda_j$ & No & None, since $|\Sigma_{1m}|\ge\lambda^2$ while $\rho_{\epsilon,C}(m-1)\to0$\\
Blocks of growing size $b$, fixed $r>0$ & $\Sigma_{ij}=r$ within a block & No & None, since $C\ge r\log^{1+\epsilon}(b-1+e)\to\infty$ is needed\\
\bottomrule
\end{tabular}
\end{table}

\paragraph{Models Covered.}
(i). \emph{AR(1) model.} Let $\xi_i=a\xi_{i-1}+\sqrt{1-a^2}\,Z_i$, $i=1,2,\ldots$, be stationary with $|a|<1$. 
Then $\rho(k)=a^k$. Since $|a|^k\log^{1+\epsilon}(k+e)\to0$, the supremum of $|\rho(k)|\log^{1+\epsilon}(k+e)$ over $k\ge1$ is finite, so $\Sigma\in\mathcal C_m(\epsilon,C)$ for every $\epsilon>0$ with a suitable $C$. For $a<0$ the correlations alternate in sign.

(ii). \emph{Moving average model.} Let $\xi_i=\sum_{j=0}^qb_jZ_{i-j}$, $i=1,2,\ldots$, be stationary with $\sum_{j=0}^qb_j^2=1$. Then $\rho(k)=\sum_{j=0}^{q-k}b_jb_{j+k}$ for $1\le k\le q$, and $\rho(k)=0$ for $k>q$. Since $|\rho(k)|\le1$, we have $\Sigma\in\mathcal C_m(\epsilon,C)$ with $C=\log^{1+\epsilon}(q+e)$
for every $\epsilon>0$.

(iii). \emph{Fractional Gaussian noise model.} Let $\xi_i=B_H(i)-B_H(i-1)$ be a stationary process, where $B_H$ is fractional Brownian motion with Hurst parameter $H\in(1/2,1)$. Then 
\[
\rho(k)=\frac12[(k+1)^{2H}-2k^{2H}+(k-1)^{2H}]\sim H(2H-1)k^{2H-2}.
\]
Since $k^{2H-2}\log^{1+\epsilon}(k+e)\to0$, we have $\Sigma\in\mathcal C_m(\epsilon,C)$ for every $\epsilon>0$ with a suitable $C$. These correlations are not summable, so Theorem~\ref{thm:main} allows long-range dependence.

(iv). \emph{Independent blocks of bounded size.} Partition $\{1,\ldots,m\}$ into consecutive blocks of size at most $B$. Let the subvectors of $\xi$ on different blocks be independent, with an arbitrary correlation matrix within each block. Then $\Sigma_{ij}=0$ for $|i-j|\ge B$, since such $i$ and $j$ lie in different blocks. For $1\le|i-j|\le B-1$, we have $|\Sigma_{ij}|\le1=\rho_{\epsilon,C}(|i-j|)$ with $C=\log^{1+\epsilon}(B-1+e)$. Hence $\Sigma\in\mathcal C_m(\epsilon,C)$. That is, the correlation matrix belongs to the class being considered.

\paragraph{Models Not Covered.}
(i). \emph{Equicorrelation.} Let $\xi_i=\sqrt r\,Z_0+\sqrt{1-r}\,Z_i$ with a fixed $r\in(0,1]$. Then $\rho(k)=r$ for all $k\ge1$. Since $r\log^{1+\epsilon}(k+e)\to\infty$, there is no $C$ with $\Sigma\in\mathcal C_m(\epsilon,C)$ for all $m$.

(ii). \emph{One-factor models.} Let $\xi_i=\lambda_iZ_0+\sqrt{1-\lambda_i^2}\,Z_i$ with $|\lambda_i|\le1$. Then $\mathrm{Corr}(\xi_i,\xi_j)=\lambda_i\lambda_j$. If $|\lambda_i|\ge\lambda>0$ for all $i$, then $|\Sigma_{1m}|\ge\lambda^2$, while $\rho_{\epsilon,C}(m-1)\to0$. So the correlation matrix is not in $\mathcal C_m(\epsilon,C)$ for large $m$, whatever $C$ is. This includes the rank-one models of \citet{LyuWang2026}, where $|\lambda_i|=1$. The common-factor models of \citet{Lei2026} have the same form, and their correlations also do not decay with the distance.

(iii). \emph{Growing blocks.} Let a block of size $b$ have a fixed within-block correlation $r>0$. Its two endpoints have correlation $r$ at lag $b-1$. A bound $|\Sigma_{ij}|\le\rho_{\epsilon,C}(|i-j|)$ therefore requires $C\ge r\log^{1+\epsilon}(b-1+e)$, which fails as $b\to\infty$.

In the end, we comment that AR(1) cannot be produced by any one-factor model. Every one-factor correlation matrix satisfies $\Sigma_{12}\Sigma_{34}=\Sigma_{13}\Sigma_{24}$, while AR(1) gives $a^2\ne a^4$ for $0<|a|<1$. In this sense, Theorem~\ref{thm:main} and the factor-model results cover different dependence structures.

\section{Experiments}\label{sec:numerics}
In this section, we evaluate our FDR bound (\ref{our:FDR:bound}) numerically, and compare it with the BY bound of \citet{BY2001} and with the empirical FDR of BH. We took $\alpha=0.1$ and $\pi_0=1$, the least favorable value. We considered two decay rates $\epsilon=0.5$ and $\epsilon=1$, and three scales $C=0.01$, $0.05$ and $0.1$. For one-sided case we used $k^*=2$ when $\epsilon=0.5$ and $k^*=1$ when $\epsilon=1$. For two-sided $p$-values only even Hermite orders appear, so we used $k^*=2$ in both cases. In every case $\epsilon k^*\ge1$.
For each $m$, we computed our FDR bound provided in the formula (\ref{our:FDR:bound})
using the envelope $\rho_{\epsilon,C}$. The constant $K$ was computed from Plackett's formula \citep{Plackett1954}. The quantities $E_{\rho_{\epsilon,C},m}(\alpha)$ were computed directly. The $\theta$ therein was selected to minimize the constant. 
We estimated the FDR of BH by Monte Carlo. We took the global null $\mu=0$ so that $\pi_0=1$. The correlation was chosen as $\rho=\rho_{\epsilon,C}$. Under the global null, the FDR equals the probability of at least one rejection. The one-sided $p$-values are $1-\Phi(T_i)$ and the two-sided ones are $2(1-\Phi(|T_i|))$. We examined $m=2^a$ for $a=2,\ldots,20$. FDR was calculated through 20000 independent replications.
 
\begin{figure}[htp]
\centering
\includegraphics[width=\textwidth]{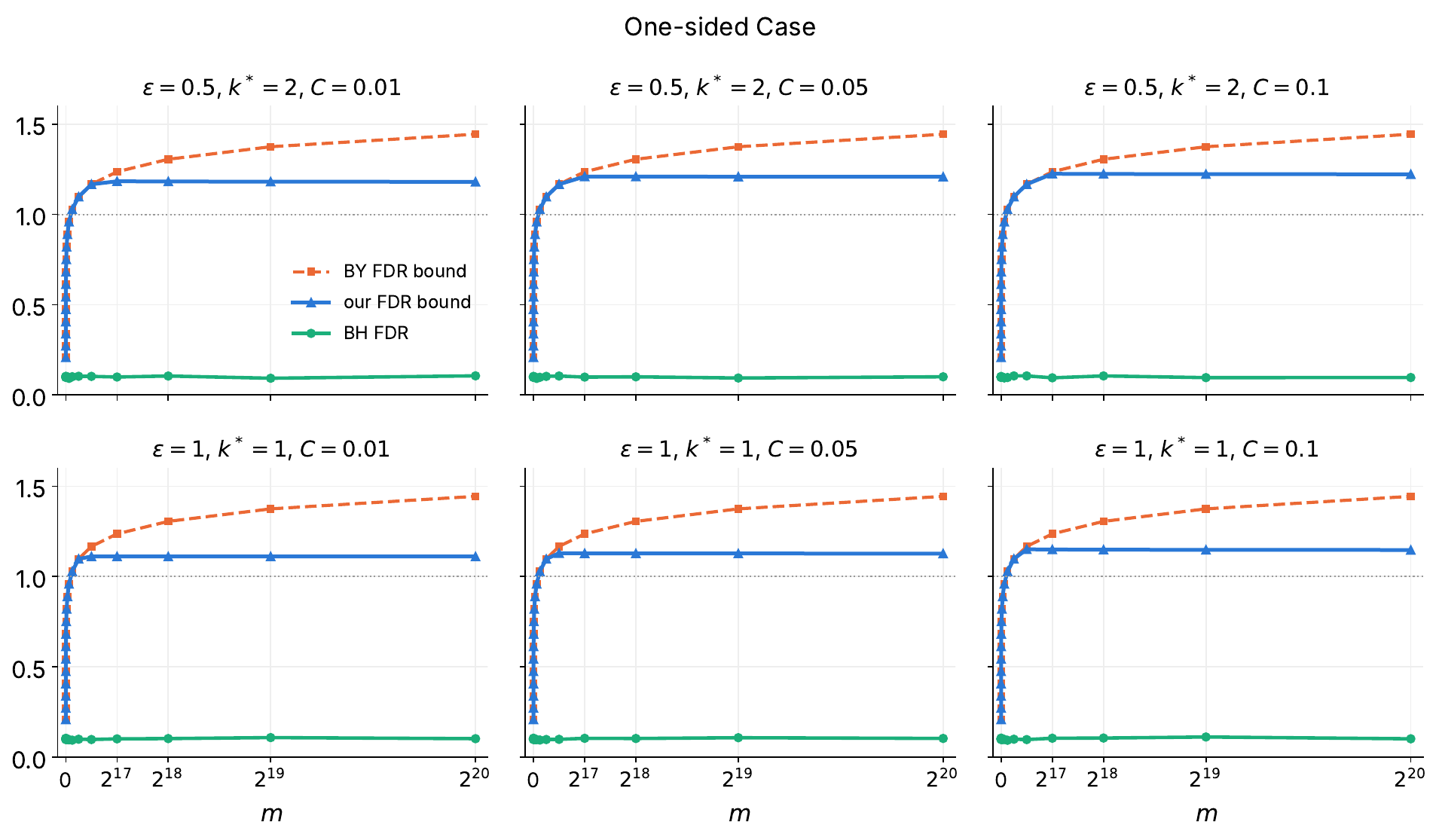}
\caption{\it One-Sided Case. FDR of BH and its upper bounds as functions of $m$. Dashed: BY FDR bound. Solid: our FDR bound. Dots: BH FDR. The dotted line marks one.}
\label{fig:one}
\end{figure}

\begin{figure}[htp]
\centering
\includegraphics[width=\textwidth]{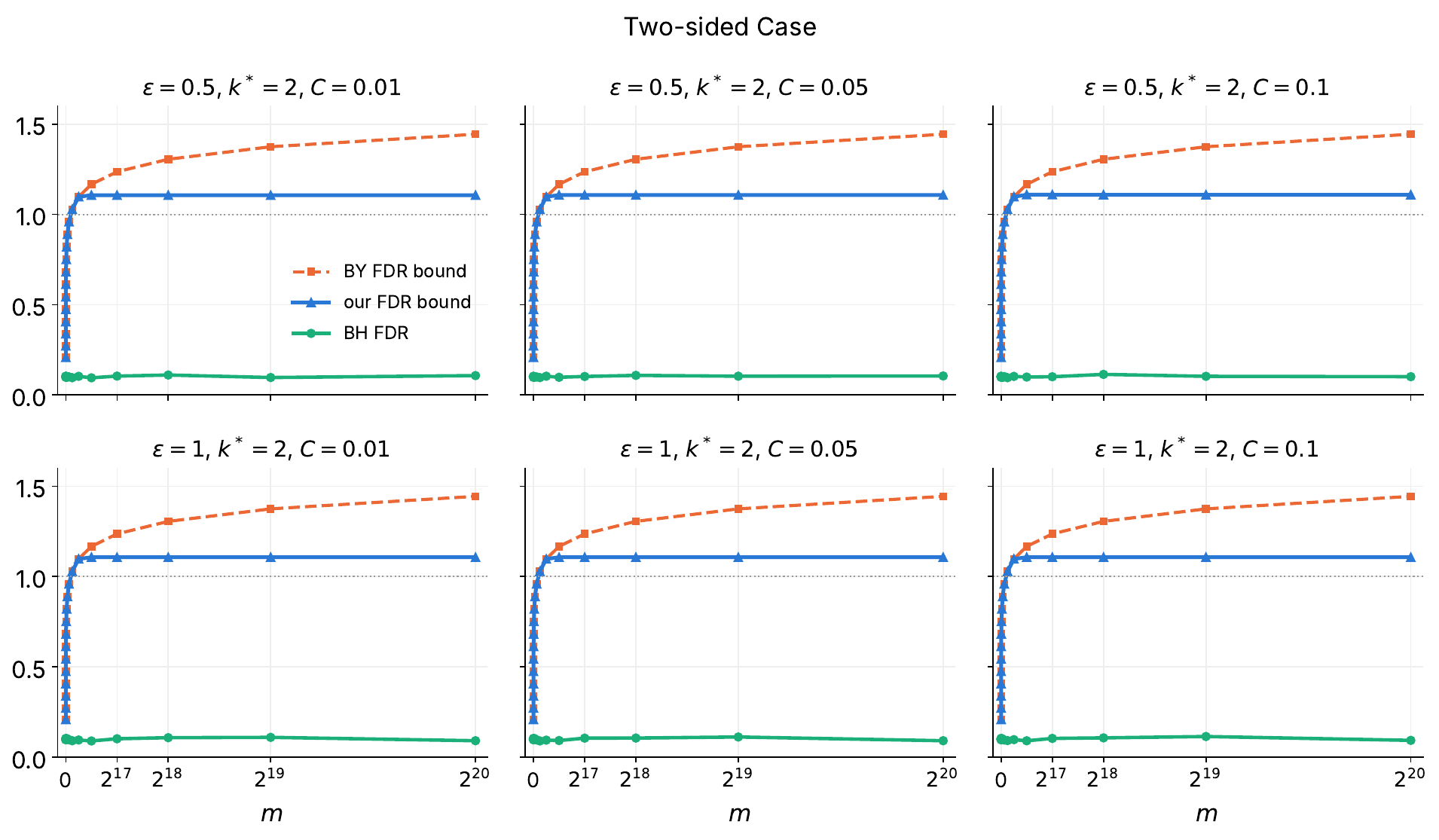}
\caption{\it Two-Sided Case. FDR of BH and its upper bounds as functions of $m$. Dashed: BY FDR bound. Solid: our FDR bound. Dots: BH FDR. The dotted line marks one.}
\label{fig:two}
\end{figure}

Figures~\ref{fig:one} and~\ref{fig:two} show the results. The BY bound of \citet{BY2001} grows like $\alpha\log m$ without limit and reaches $1.44$ at $m=2^{20}$. Our bound coincides with it only for small $m$ and stays flat for large $m$. In the one-sided case it separates from it between $m=2^{13}$ and $2^{14}$ for both $\epsilon=0.5$ and $\epsilon=1$. After that it stays flat, at $0.97$ to $0.98$ when $\epsilon=1$. When $\epsilon=0.5$ it stays at $0.996$ for $C=0.01$ and reaches the trivial value $1$ for $C=0.05$ and $0.1$. A larger $C$ or a smaller $\epsilon$ means stronger correlation, and it moves the bound up and the separation to the right. In the two-sided case all six bounds are nearly the same. They separate from the BY bound between $m=2^{13}$ and $2^{14}$ and stay at $0.97$. The empirical FDR lies between $0.09$ and $0.11$ in all cases and shows no growth in $m$.
 
\paragraph{Acknowledgment.} Claude Opus 5.5 (Anthropic) assisted in developing several of the proofs, such as the Hermite expansion and covariance bounds (Lemmas~\ref{lem:hermite}--\ref{lem:A}), and in editing the manuscript. All results were checked by the authors, who takes full responsibility for the content.

\section{Appendix: Proofs}\label{sec:appendix}

This section proves Theorem~\ref{thm:main}. We use the notation of Section~\ref{sec:main} throughout. The argument first works with an arbitrary correlation matrix through its lag envelope $\rho$, defined below, and then specializes to the envelope $\rho_{\epsilon,C}$. Section~\ref{sec:prelim} extends the quantities of Section~\ref{sec:main} to an arbitrary $\rho$ and reduces the problem to a centered Gaussian sequence. Section~\ref{sec:gauss} collects Gaussian and Hermite facts. Section~\ref{sec:var} bounds a variance. Section~\ref{sec:masterineq} combines these into an FDR bound for a single correlation matrix (Proposition~\ref{prop:master}). Section~\ref{sec:proofmain} evaluates this bound on the envelope. In Sections~\ref{sec:prelim}--\ref{sec:masterineq}, $k^*\in\mathcal S$ is arbitrary. The condition $\epsilon k^*\ge1$ is used only in Section~\ref{sec:proofmain}.

\subsection{Notation and Reduction to a Centered Sequence}\label{sec:prelim}
\emph{Lag envelope of a correlation matrix.} For an $m\times m$ correlation matrix $\Sigma$, define its lag envelope $\rho_\Sigma\colon\mathbb N_0\to[0,1]$ by $\rho_\Sigma(0)=1$, $\rho_\Sigma(l)=\max_{|i-j|=l}|\Sigma_{ij}|$ for $1\le l\le m-1$, and $\rho_\Sigma(l)=0$ for $l\ge m$. Then $|\Sigma_{ij}|\le\rho_\Sigma(|i-j|)$ for all $i\ne j$. Moreover, $\Sigma\in\mathcal C_m(\epsilon,C)$ if and only if $\Sigma\succeq0$, $\Sigma_{ii}=1$ and $\rho_\Sigma\le\rho_{\epsilon,C}$ pointwise.

\emph{General correlation sequences.} For $\rho\colon\mathbb N_0\to[-1,1]$ with $\rho(0)=1$, let $D_\rho=\#\{l\ge1:|\rho(l)|>1/2\}$, and let $K_\rho$, $G_{\rho,k}(m)$, $e_{\rho,m}(z)$ and $\Lambda_{\rho,m}(\alpha)$ be defined as $K_{\epsilon,C}$, $G_{\epsilon,C,k}(m)$, $e_{\epsilon,C,m}(z)$ and $\Lambda_{\epsilon,C,m}(\alpha)$ in Section~\ref{sec:main}, with $\rho_{\epsilon,C}(l)$ replaced by $|\rho(l)|$ for every $l\ge1$. That is,
\[
K_\rho=1+\sup_{z\ge0}\sum_{l\ge1:\,|\rho(l)|>r_z}\tau(\rho(l),z),\qquad G_{\rho,k}(m)=\sum_{l=1}^{m-1}|\rho(l)|^k,
\]
and $e_{\rho,m}$ and $\Lambda_{\rho,m}(\alpha)$ are obtained from $e_{\epsilon,C,m}$ and $\Lambda_{\epsilon,C,m}(\alpha)$ by replacing $G_{\epsilon,C,k}(m)$ with $G_{\rho,k}(m)$. 
Specifically,  for $z\ge0$,
\begin{eqnarray*}
e_{\rho,m}(z)&=&2\Big(\frac{\phi(z)}{1-\Phi(z)}\Big)^2\mathsf A(z)\,\frac{G_{\rho,k^*}(m)}m,\\
\Lambda_{\rho,m}(\alpha)&=&\left\{\begin{array}{ll}\Big(\sum_{k\in\mathcal S,\ k<k^*}b_k\Big[k!\Big(\frac1m+\frac{2G_{\rho,k}(m)}m\Big)\Big]^{1/2}\Big)^2, & k^*>\min\mathcal S,\\
0, & k^*=\min\mathcal S.\end{array}\right.
\end{eqnarray*}
Since $\rho_{\epsilon,C}\ge0$, we have $K_{\rho_{\epsilon,C}}=K_{\epsilon,C}$, $e_{\rho_{\epsilon,C},m}=e_{\epsilon,C,m}$ and $\Lambda_{\rho_{\epsilon,C},m}(\alpha)=\Lambda_{\epsilon,C,m}(\alpha)$. For $z\le\sqrt2$ we have $r_z=1/2$, so the sum in $K_\rho$ runs over the $D_\rho$ lags with $|\rho(l)|>1/2$ and equals $2D_\rho$. Hence $K_\rho\ge1+2D_\rho\ge1$.

\emph{Reduction.}
Throughout, $\Sigma$ is an $m\times m$ correlation matrix, $\rho=\rho_\Sigma$, and $\xi=T-\mu\sim N(0,\Sigma)$. Let $p(x)=1-\Phi(x)$ (one-sided) or $p(x)=2(1-\Phi(|x|))$ (two-sided), so that $P_i=p(\xi_i)$ for $i\in\mathcal H_0$. For $w\in(0,1/2]$ we have $z_w\ge0$, and $p(\xi_i)\le w$ if and only if $\xi_i\ge z_w$ (one-sided), respectively $|\xi_i|\ge z_w$ (two-sided). Put
\[
A_i(w)=\mathbf 1\{p(\xi_i)\le w\},\qquad \Ntil(w)=\sum_{i\in\mathcal H_0}A_i(w),\qquad \Stil=\sup_{u\in[\alpha/m,\alpha]}\frac{\Ntil(u)}{mu},
\]
so that $\E A_i(w)=w$ and $\E\Ntil(w)=m_0w=\pi_0mw$. Note that $\Ntil$ counts only the true nulls, and $\Stil$ depends on $\mu$ only through $\mathcal H_0$. If $m_0=0$, then $V=0$ and $\FDR=0$, so we assume $m_0\ge1$ whenever we use $\Stil$.

\begin{lemma}\label{lem:reduction}
$\FDP\le\min\{1,\alpha\Stil\}$ almost surely. Consequently $\FDR\le\alpha\int_0^{1/\alpha}\P(\Stil\ge t)\,dt$.
\end{lemma}
\begin{proof}
If $R=0$ then $\FDP=0$. Otherwise $u_R=\alpha R/m\in[\alpha/m,\alpha]$, and every rejected true null satisfies $p(\xi_i)=P_i\le u_R$, so $V\le\Ntil(u_R)$ and $\FDP=V/R=\alpha V/(mu_R)\le\alpha\Stil$. Finally $\E\min\{1,Y\}=\int_0^1\P(Y\ge s)\,ds$ with $Y=\alpha\Stil$.
Proof is completed.
\end{proof}

Recall $J=\lceil\log_2m\rceil$ and $w_j=2^{-j}\alpha$, and let also $w_J=2^{-J}\alpha$ and $a_j=mw_j$ for $0\le j\le J$. All levels satisfy $w_j\le\alpha\le1/2$.

\begin{lemma}\label{lem:dyadic}
For every $t>0$, $\{\Stil\ge t\}\subseteq\bigcup_{j=0}^{J-1}\{\Ntil(w_j)\ge\frac t2a_j\}$.
\end{lemma}
\begin{proof}
Since $w_J\le\alpha/m$, every $u\in[\alpha/m,\alpha]$ lies in some $[w_{j+1},w_j]$ with $0\le j\le J-1$, and then $\Ntil(u)/(mu)\le\Ntil(w_j)/(mw_{j+1})=2\Ntil(w_j)/a_j$. Proof is completed.
\end{proof}

\begin{lemma}\label{lem:geom}
If $a_0>\dots>a_{J-1}>0$ with $a_j=2a_{j+1}$ and $b>0$, then $\sum_{j=0}^{J-1}\min\{a_j,b/a_j\}\le3\sqrt b$.
\end{lemma}
\begin{proof}
The terms with $a_j\le\sqrt b$ form a geometric sequence with ratio $1/2$ and largest term $a_*\le\sqrt b$, so they sum to at most $2a_*$. The terms with $a_j>\sqrt b$ are at most $b/a_j$, a geometric sequence with ratio $1/2$ and largest term $b/a^*$, where $a^*>\sqrt b$ is the smallest such $a_j$; they sum to at most $2b/a^*<2\sqrt b$. If only one group is present, the total is at most $2\sqrt b$. If both are present, $a^*=2a_*$, the total is at most $2a_*+b/a_*$, and $a_*\in[\sqrt b/2,\sqrt b]$; the convex function $x\mapsto2x+b/x$ equals $3\sqrt b$ at both endpoints. Proof is completed.
\end{proof}

\subsection{Gaussian and Hermite Facts}\label{sec:gauss}
Recall $He_k$, $R(r,z)$, $r_z$ and $\mathsf A(z)$ from Section~\ref{sec:main}. Let $c_k(z)=\phi(z)He_{k-1}(z)/k!$ for $k\ge1$.
 
\begin{lemma}\label{lem:hermite}
Let $z\ge0$ and $|r|<1$.
\begin{enumerate}[label=(\roman*).,ref=(\roman*)]
\item\label{lem:hermite:i} With $X\sim N(0,1)$, 
\begin{eqnarray*}
\mathbf 1\{X>z\}&=&1-\Phi(z)+\sum_{k\ge1}c_k(z)He_k(X),\\ 
\mathbf 1\{|X|>z\}&=&2(1-\Phi(z))+\sum_{k\in2\mathbb N}2c_k(z)He_k(X),
\end{eqnarray*}
where $c_k(z)=\frac{\phi(z)He_{k-1}(z)}{k!}$ and the series converge in mean square. For example,
\[
\textrm{$\E\Big[\Big(\mathbf 1\{X>z\}-\big(1-\Phi(z)\big)-\sum_{k=1}^nc_k(z)He_k(X)\Big)^2\Big]\to0$ as $n\to\infty$.}
\]
\item\label{lem:hermite:ii} If $(X,Y)$ is standard bivariate normal with correlation $r$, then for any $k,l\ge1$, 
\[
\E He_k(X)He_l(Y)=\delta_{kl}\,k!\,r^k.
\]
Hence
\begin{eqnarray*}
R(r,z)&=&\phi(z)^2\sum_{k\ge1}\frac{He_{k-1}(z)^2}{k!}r^k,\\
\Cov(\mathbf 1\{|X|>z\},\mathbf 1\{|Y|>z\})&=&4\phi(z)^2\sum_{k\in2\mathbb N}\frac{He_{k-1}(z)^2}{k!}r^k .
\end{eqnarray*}
\item\label{lem:hermite:iii} $R(r,z)=\phi(z)^2\int_0^r e^{z^2s/(1+s)}(1-s^2)^{-1/2}\,ds$. Consequently $R(r,z)$ has the sign of $r$, $R(\cdot,z)$ is increasing, $|R(r,z)|\le R(|r|,z)$, and the two-sided covariance in \ref{lem:hermite:ii} equals $2[R(r,z)+R(-r,z)]\in[0,2R(|r|,z)]$.
\end{enumerate}
\end{lemma}
\begin{proof}
\emph{\ref{lem:hermite:i}}. Differentiating $\phi^{(k-1)}=(-1)^{k-1}He_{k-1}\phi$ gives
\[
\big(He_{k-1}\phi\big)'=(-1)^{k-1}\phi^{(k)}=-He_k\phi .
\]
Since $He_{k-1}(x)\phi(x)\to0$ as $x\to\infty$, for $k\ge1$,
\[
\E\big[\mathbf 1\{X>z\}He_k(X)\big]=\int_z^\infty He_k(x)\phi(x)\,dx=-\Big[He_{k-1}(x)\phi(x)\Big]_z^\infty=\phi(z)He_{k-1}(z),
\]
while $\E[\mathbf 1\{X>z\}He_0(X)]=1-\Phi(z)$. The polynomials $He_k/\sqrt{k!}$, $k\ge0$, form an orthonormal basis of the space of functions $f$ with $\E f(X)^2<\infty$, equipped with the inner product $\E[f(X)g(X)]$; in particular $\E He_k(X)^2=k!$. Since $\mathbf 1\{x>z\}$ is bounded, it equals in mean square its expansion in this basis, and the coefficient of $He_k$ is
\[
\frac{\E[\mathbf 1\{X>z\}He_k(X)]}{\E He_k(X)^2}=\frac{\phi(z)He_{k-1}(z)}{k!}=c_k(z),\qquad k\ge1.
\]
This proves the first expansion. For the second, $\mathbf 1\{|x|>z\}=\mathbf 1\{x>z\}+\mathbf 1\{-x>z\}$ because $z\ge0$. Since $X$ and $-X$ have the same distribution and $He_k(-x)=(-1)^kHe_k(x)$,
\[
\E\big[\mathbf 1\{-X>z\}He_k(X)\big]=\E\big[\mathbf 1\{X>z\}He_k(-X)\big]=(-1)^k\phi(z)He_{k-1}(z).
\]
Hence the coefficient of $He_k$ in the expansion of $\mathbf 1\{|x|>z\}$ is $(1+(-1)^k)c_k(z)$, which equals $2c_k(z)$ for even $k$ and $0$ for odd $k$, and the constant term is $2(1-\Phi(z))$.
 
\noindent\emph{\ref{lem:hermite:ii}}. The generating function of the Hermite polynomials is $e^{sx-s^2/2}=\sum_{k\ge0}He_k(x)s^k/k!$. Write $Y=rX+\sqrt{1-r^2}\,W$ with $W\sim N(0,1)$ independent of $X$, 
which does not change the joint distribution of $(X,Y)$.
For real $s,t$,
\[
\begin{aligned}
\E\big[e^{sX-s^2/2}e^{tY-t^2/2}\big]&=e^{-(s^2+t^2)/2}\,\E e^{(s+rt)X}\,\E e^{t\sqrt{1-r^2}\,W}\\
&=e^{-(s^2+t^2)/2}e^{(s+rt)^2/2}e^{t^2(1-r^2)/2}=e^{rst}.
\end{aligned}
\]
Expanding both sides in powers of $s$ and $t$ and comparing the coefficients of $s^kt^l$ gives
\[
\frac{\E[He_k(X)He_l(Y)]}{k!\,l!}=\delta_{kl}\frac{r^k}{k!},
\]
that is,
\begin{equation}\label{eq:ONB}
\E He_k(X)He_l(Y)=\delta_{kl}\,k!\,r^k.
\end{equation}
 
Next let $F_n(x)=\sum_{k=1}^nc_k(z)He_k(x)$. By \ref{lem:hermite:i}, $F_n(X)\to\mathbf 1\{X>z\}-(1-\Phi(z))$ and $F_n(Y)\to\mathbf 1\{Y>z\}-(1-\Phi(z))$ in mean square, so by the Cauchy--Schwarz inequality $\E[F_n(X)F_n(Y)]\to R(r,z)$. Therefore
\[
R(r,z)=\lim_{n\to\infty}\sum_{k,l=1}^nc_k(z)c_l(z)\,\E\big[He_k(X)He_l(Y)\big]=\sum_{k\ge1}k!\,c_k(z)^2r^k=\phi(z)^2\sum_{k\ge1}\frac{He_{k-1}(z)^2}{k!}r^k.
\]
The two-sided formula follows in the same way from the second expansion in \ref{lem:hermite:i}, with $c_k$ replaced by $2c_k$ for even $k$ and by $0$ for odd $k$.
 
\noindent\emph{\ref{lem:hermite:iii}}. Let $\phi_2(x,y;r)$ denote the standard bivariate normal density with correlation $r$. It satisfies $\partial_r\phi_2=\partial_x\partial_y\phi_2$ \citep{Plackett1954}, so
\[
\frac{\partial}{\partial r}\P(X>z,Y>z)=\int_z^\infty\!\!\int_z^\infty\partial_x\partial_y\phi_2(x,y;r)\,dx\,dy=\phi_2(z,z;r)=\frac{1}{2\pi\sqrt{1-r^2}}\exp\Big\{-\frac{z^2}{1+r}\Big\}.
\]
Since $-z^2/(1+r)=-z^2+z^2r/(1+r)$ and $\phi(z)^2=e^{-z^2}/(2\pi)$, the right-hand side equals $\phi(z)^2e^{z^2r/(1+r)}(1-r^2)^{-1/2}$. At $r=0$, $\P(X>z,Y>z)=(1-\Phi(z))^2$, that is, $R(0,z)=0$. Integrating from $0$ to $r$ gives the formula for $R(r,z)$. The integrand is positive, so $R(r,z)$ has the sign of $r$ and $R(\cdot,z)$ is increasing. For $r<0$, substituting $s\mapsto-s$,
\[
|R(r,z)|=\phi(z)^2\int_0^{|r|}\frac{e^{-z^2s/(1-s)}}{\sqrt{1-s^2}}\,ds\le\phi(z)^2\int_0^{|r|}\frac{e^{z^2s/(1+s)}}{\sqrt{1-s^2}}\,ds=R(|r|,z).
\]
For the two-sided covariance, write $\mathbf 1\{|X|>z\}=\mathbf 1\{X>z\}+\mathbf 1\{-X>z\}$ and similarly for $Y$, and expand the covariance into four terms. The pairs $(X,Y)$ and $(-X,-Y)$ have correlation $r$, and the pairs $(X,-Y)$ and $(-X,Y)$ have correlation $-r$, so
\[
\Cov\big(\mathbf 1\{|X|>z\},\mathbf 1\{|Y|>z\}\big)=2\big[R(r,z)+R(-r,z)\big].
\]
By \ref{lem:hermite:ii} this is a series in even powers of $r$ with nonnegative coefficients, hence nonnegative. Finally, $R(r,z)+R(-r,z)=R(|r|,z)+R(-|r|,z)\le R(|r|,z)$ because $R(-|r|,z)\le0$. Proof is completed.
\end{proof}

\begin{lemma}\label{lem:mills}
Let $z\ge0$ and $h(z)=(z+\sqrt{z^2+4})/2$.
\begin{enumerate}[label=(\roman*).,ref=(\roman*)]
\item\label{lem:mills:i} $\phi(z)/(1-\Phi(z))\le h(z)$ \citep{Birnbaum1942}, $h(z)^2\le z^2+2$, $h(z)\le1+z$, and $1-\Phi(z)\le\frac12e^{-z^2/2}$. Hence $z_w^2\le2\log\frac1{2w}$ (one-sided) and $z_w^2\le2\log\frac1w$ (two-sided).
\item\label{lem:mills:ii} If $|r|\le1/2$, then $R(|r|,z)/(1-\Phi(z))\le(1+z)e^{-z^2/6}/\sqrt{6\pi}$.
\item\label{lem:mills:iii} $|He_n(z)|\le(z+\sqrt{n+1})^n$ for $n\ge0$.
\end{enumerate}
\end{lemma}
\begin{proof}
\ref{lem:mills:i}. $h^2=\frac12(z^2+2+z\sqrt{z^2+4})\le z^2+2$ since $z^2(z^2+4)\le(z^2+2)^2$. Others are standard. 

\noindent\ref{lem:mills:ii}. By Lemma~\ref{lem:hermite}\ref{lem:hermite:iii}, $s/(1+s)\le1/3$ and $(1-s^2)^{-1/2}\le2/\sqrt3$ for $0\le s\le1/2$, so $R(|r|,z)\le\phi(z)^2e^{z^2/3}/\sqrt3$; combine with $1-\Phi(z)\ge\phi(z)/(1+z)$ from \ref{lem:mills:i}. 

\noindent\ref{lem:mills:iii}. $He_n(z)=\E(z+\mathrm iY)^n$ with $Y\sim N(0,1)$, so by Minkowski's inequality $|He_n(z)|\le(z+\|Y\|_n)^n$; moreover $\|Y\|_q^q=(q-1)!!\le q^{q/2}$ for even $q$, so $\|Y\|_n\le\sqrt{n+1}$ (use $\|Y\|_n\le\|Y\|_{n+1}$ for odd $n$). Proof is completed.
\end{proof}

By the definition of $\mathsf A(z)$, termwise, for $|r|\le r_z$,
\begin{equation}\label{eq:remainder}
\sum_{k\in\mathcal S,\ k\ge k^*}\frac{He_{k-1}(z)^2}{k!}|r|^k\le\mathsf A(z)\,|r|^{k^*}.
\end{equation}
The following lemma provides an upper bound for $\mathsf{A}(z)$ for various choices of $k^*$.
\begin{lemma}\label{lem:A}
Let $z\ge0$.
\begin{enumerate}[label=(\roman*).,ref=(\roman*)]
\item\label{lem:A:i} In one-sided case, if we choose $k^*=1$ (the smallest element of $\mathcal S$), then $\mathsf A(z)\le a_1:=2(e-1)/\sqrt3<1.985$.
\item\label{lem:A:ii} In two-sided case, if we choose $k^*=2$ (the smallest element of $\mathcal S$), then $\mathsf A(z)\le a_2z^2$ with $a_2:=2/\sqrt3<1.155$.
\item\label{lem:A:iii} For every $k^*\in\mathcal S$, in either one- or two-sided case, we have $\mathsf A(z)\le a_1\max\{2,z^2\}^{k^*-1}$.
\end{enumerate}
\end{lemma}
In \ref{lem:A:i} and \ref{lem:A:ii} the series defining $\mathsf A(z)$ is a full covariance series, which Plackett's formula sums in closed form. Part \ref{lem:A:iii} compares $\mathsf A(z)$ with the same series, and it reduces to \ref{lem:A:i} when $k^*=1$.
\begin{proof}
Write $x=z^2r_z\le1$ and use $(1-s^2)^{-1/2}\le2/\sqrt3$ for $0\le s\le r_z\le1/2$.

\noindent \ref{lem:A:i}. By Lemma~\ref{lem:hermite}\ref{lem:hermite:ii},\ref{lem:hermite:iii}, as well as $(e^x-1)/x$ is increasing, it holds that
\[
\mathsf A(z)=R(r_z,z)/(\phi^2r_z)\le\frac{2}{\sqrt3}\,r_z^{-1}\int_0^{r_z}e^{z^2s}ds=\frac2{\sqrt3}\,\frac{e^x-1}{x}\le\frac{2(e-1)}{\sqrt3}.
\]

\noindent \ref{lem:A:ii}. By Lemma~\ref{lem:hermite}, $\mathsf A(z)=[R(r_z,z)+R(-r_z,z)]/(2\phi^2r_z^2)$ and
\[
R(r,z)+R(-r,z)=\phi^2\int_0^r(1-s^2)^{-1/2}\big(e^{z^2s/(1+s)}-e^{-z^2s/(1-s)}\big)ds.
\]
Using $e^u-e^{-v}\le ue^u+v$ with $u=z^2s/(1+s)\le z^2s$ and $v=z^2s/(1-s)\le2z^2s$, we have
\[
\mathsf A(z)\le\frac{1}{\sqrt3\,r_z^2}\int_0^{r_z}z^2s\,(e^{z^2s}+2)\,ds=\frac{z^2}{\sqrt3}\Big(\frac{xe^x-e^x+1}{x^2}+1\Big)\le\frac{2z^2}{\sqrt3},
\]
because $(xe^x-e^x+1)/x^2$ is increasing and equals $1$ at $x=1$.

\noindent \ref{lem:A:iii}. Fix $k^*\in\mathcal S$. Since $r_z^{\,k-k^*}=r_z^{-k^*}r_z^k$ and all terms below are nonnegative,
\[
\mathsf A(z)=r_z^{-k^*}\sum_{k\in\mathcal S,\ k\ge k^*}\frac{He_{k-1}(z)^2}{k!}\,r_z^{\,k}\le r_z^{-k^*}\sum_{k\ge1}\frac{He_{k-1}(z)^2}{k!}\,r_z^{\,k}=r_z^{-k^*}\frac{R(r_z,z)}{\phi(z)^2},
\]
where the last step is Lemma~\ref{lem:hermite}\ref{lem:hermite:ii}. By Lemma~\ref{lem:hermite}\ref{lem:hermite:iii}, as in \ref{lem:A:i},
\[
\frac{R(r_z,z)}{\phi(z)^2}\le\frac2{\sqrt3}\int_0^{r_z}e^{z^2s}\,ds=\frac2{\sqrt3}\,r_z\,\frac{e^x-1}x\le\frac{2(e-1)}{\sqrt3}\,r_z=a_1r_z.
\]
Hence $\mathsf A(z)\le a_1r_z^{1-k^*}=a_1\max\{2,z^2\}^{k^*-1}$, because $r_z^{-1}=\max\{2,z^2\}$. Proof is completed.
\end{proof}

\subsection{A Variance Bound}\label{sec:var}
Fix $k^*\in\mathcal S$ and a level $w\in(0,1/2]$ with threshold $z=z_w$. Let $\gamma_k=c_k(z)$ (one-sided) or $\gamma_k=2c_k(z)$ (two-sided), and
\[
Q_i(w)=A_i(w)-w-\sum_{k\in\mathcal S,\ k<k^*}\gamma_kHe_k(\xi_i).
\]
For $k^*=\min\mathcal S$ the sum is empty and $Q_i(w)=A_i(w)-w$.

\begin{proposition}\label{prop:var}
Let $z=z_w$ and let $I\subseteq\{1,\dots,m\}$ have $|I|=m_0\ge1$. Then $\Var\big(\sum_{i\in I}Q_i(w)\big)\le K_\rho\,m_0w+e_{\rho,m}(z_w)\,m_0mw^2$.
\end{proposition}
\begin{proof}
Recall that $w\le1/2$, so $z=z_w\ge0$ and Lemmas~\ref{lem:hermite}--\ref{lem:A} apply. Expanding the variance of a sum,
\begin{equation}\label{eq:varsplit}
\Var\Big(\sum_{i\in I}Q_i(w)\Big)=\sum_{i\in I}\Var Q_i(w)+\sum_{i,i'\in I,\ i\ne i'}\Cov\big(Q_i(w),Q_{i'}(w)\big).
\end{equation}
We bound the diagonal terms in Step~1, derive three bounds on a single covariance in Step~2, and sum the covariances lag by lag in Step~3.

\emph{Step 1: diagonal terms.} 
Since $He_0\equiv 1$ and (\ref{eq:ONB}), we have $\E H_k(\xi_i)=0$ for $k\in\mathcal{S}$, hence $\E Q_i(w)=0$. By Lemma \ref{lem:hermite}\ref{lem:hermite:i}, we have
\[
\Var Q_i(w)=\E Q_i(w)^2\le\E\big(A_i(w)-w\big)^2=w(1-w)\le w.
\]
The diagonal terms in \eqref{eq:varsplit} therefore contribute at most $m_0w$.

\emph{Step 2: a single covariance.} Fix $i\ne i'$ and let $r=\Sigma_{ii'}$, so that $|r|\le\rho(|i-i'|)$. Set
\[
\sigma(r)=\sum_{k\in\mathcal S,\ k\ge k^*}\frac{He_{k-1}(z)^2}{k!}\,r^k .
\]
We claim that, in both the one-sided and the two-sided case,
\begin{equation}\label{eq:covsigma}
\Cov\big(Q_i(w),Q_{i'}(w)\big)=w^2\Big(\frac{\phi(z)}{1-\Phi(z)}\Big)^2\sigma(r)\qquad\text{if }|r|<1.
\end{equation}
Indeed, by Lemma~\ref{lem:hermite}\ref{lem:hermite:i}, $Q_i(w)=\sum_{k\in\mathcal S,\,k\ge k^*}\gamma_kHe_k(\xi_i)$ and similarly for $i'$. The pair $(\xi_i,\xi_{i'})$ is standard bivariate normal with correlation $r$, so Lemma~\ref{lem:hermite}\ref{lem:hermite:ii} gives $\Cov(Q_i,Q_{i'})=\sum_{k\in\mathcal S,\,k\ge k^*}k!\,\gamma_k^2r^k$. In the one-sided case $k!\,\gamma_k^2=\phi(z)^2He_{k-1}(z)^2/k!$ and $w=1-\Phi(z)$. In the two-sided case $k!\,\gamma_k^2=4\phi(z)^2He_{k-1}(z)^2/k!$ and $w=2(1-\Phi(z))$. In both cases this is \eqref{eq:covsigma}. We now derive three bounds.

\noindent (a). \emph{Trivial bound.} By the Cauchy--Schwarz inequality and Step~1,
\[
|\Cov(Q_i,Q_{i'})|\le(\Var Q_i\,\Var Q_{i'})^{1/2}\le w\qquad\text{for every }r\in[-1,1].
\]

\noindent (b). \emph{Comparison with the full covariance.} Let $|r|\le1/2$. Dropping the restriction $k\ge k^*$ and by Lemma~\ref{lem:hermite}\ref{lem:hermite:ii},\ref{lem:hermite:iii}, we have
\[
|\sigma(r)|\le\sum_{k\in\mathcal S}\frac{He_{k-1}(z)^2}{k!}|r|^k=
\begin{cases}
R(|r|,z)/\phi(z)^2&\text{(one-sided)},\\
[R(|r|,z)+R(-|r|,z)]/(2\phi(z)^2)\le R(|r|,z)/(2\phi(z)^2)&\text{(two-sided)}.
\end{cases}
\]
Inserting this into \eqref{eq:covsigma} with $w=1-\Phi(z)$, respectively $w=2(1-\Phi(z))$, gives in both cases
\[
|\Cov(Q_i,Q_{i'})|\le w\,\frac{R(|r|,z)}{1-\Phi(z)}.
\]
Together with (a),
\begin{eqnarray*}
|\Cov(Q_i,Q_{i'})|&\le& w\min\Big\{1,\frac{R(|r|,z)}{1-\Phi(z)}\Big\},\quad\text{if }|r|\le\tfrac12,\\
|\Cov(Q_i,Q_{i'})|&\le& w,\quad\text{if }|r|>\tfrac12.
\end{eqnarray*}
By the definition of $\tau$, both cases read
\begin{equation}\label{eq:covtau}
|\Cov(Q_i,Q_{i'})|\le\frac w2\,\tau(r,z)\qquad\text{for all }r\in[-1,1].
\end{equation}

\noindent (c). 
\emph{Weak correlations.} Let $|r|\le r_z$. Then $|\sigma(r)|\le\mathsf A(z)|r|^{k^*}$ by \eqref{eq:remainder}, and \eqref{eq:covsigma} gives
\begin{equation}\label{eq:covweak}
|\Cov(Q_i,Q_{i'})|\le w^2\Big(\frac{\phi(z)}{1-\Phi(z)}\Big)^2\mathsf A(z)\,|r|^{k^*}.
\end{equation}
Bound (c) is quadratic in $w$ and of order $|r|^{k^*}$, so it is effective for the many weakly correlated pairs. Bound (b) is linear in $w$ and is used for the few remaining pairs.

\emph{Step 3: summing over pairs.} Group the off-diagonal pairs by their lag $l=|i-i'|\in\{1,\dots,m-1\}$. Each $i\in I$ has at most two partners $i'$ with $|i-i'|=l$, so for each $l$ there are at most $2m_0$ ordered pairs $(i,i')$ in $I$ with lag $l$, all with absolute correlation at most $\rho(l)$. Split the lags into
\[
L_1=\{l:\rho(l)>r_z\},\qquad L_2=\{l\le m-1:\rho(l)\le r_z\}.
\]
For $l\in L_1$ we use \eqref{eq:covtau} and $\tau(\Sigma_{ii'},z)\le\tau(\rho(l),z)$, which holds because $\tau(r,z)$ is nondecreasing in $|r|$ by Lemma~\ref{lem:hermite}\ref{lem:hermite:iii}:
\[
\sum_{l\in L_1}\ \sum_{i,i'\in I,\,|i-i'|=l}|\Cov(Q_i,Q_{i'})|\le\sum_{l\in L_1}2m_0\cdot\frac w2\,\tau(\rho(l),z)=m_0w\,S_\rho(z)\le m_0w\,(K_\rho-1),
\]
where $S_\rho(z)=\sum_{l:|\rho(l)|>r_z}\tau(\rho(l),z)$ and $K_\rho=1+\sup_{z'\ge0}S_\rho(z')$. For $l\in L_2$, every pair with lag $l$ has $|\Sigma_{ii'}|\le\rho(l)\le r_z$, so we use \eqref{eq:covweak} and $|\Sigma_{ii'}|^{k^*}\le\rho(l)^{k^*}$:
\[
\sum_{l\in L_2}\ \sum_{i,i'\in I,\,|i-i'|=l}|\Cov(Q_i,Q_{i'})|\le2m_0\,w^2\Big(\frac{\phi(z)}{1-\Phi(z)}\Big)^2\mathsf A(z)\sum_{l=1}^{m-1}\rho(l)^{k^*}=e_{\rho,m}(z_w)\,m_0mw^2,
\]
since $\sum_{l=1}^{m-1}\rho(l)^{k^*}=G_{\rho,k^*}(m)$ and $2m_0w^2G_{\rho,k^*}(m)=m_0mw^2\cdot2G_{\rho,k^*}(m)/m$.

Inserting Steps~1 and~3 into \eqref{eq:varsplit} gives
\[
\Var\Big(\sum_{i\in I}Q_i(w)\Big)\le m_0w+(K_\rho-1)\,m_0w+e_{\rho,m}(z_w)\,m_0mw^2=K_\rho\,m_0w+e_{\rho,m}(z_w)\,m_0mw^2.
\]
Proof is completed.
\end{proof}

\subsection{A Useful Inequality}\label{sec:masterineq}
The following lemma collects two properties of $\Psi_{\pi_0}$ in \eqref{eq:def:Psi}. We write $g(t)=g(t;K,E,\pi_0)$.
\begin{lemma}\label{lem:Psi}
Let $K\ge1$, $E\ge0$ and $\alpha\in(0,1/2]$.
\begin{enumerate}[label=(\roman*).,ref=(\roman*)]
\item\label{lem:Psi:ii} $\Psi_{\pi_0}(K,E,\alpha)$ is nondecreasing in $K$, in $E$ and in $\pi_0$, and $2\pi_0\le\Psi_{\pi_0}(K,E,\alpha)\le1/\alpha$.
\item\label{lem:Psi:iii} Let $\pi_0\in(0,1]$. For every $s>0$,
\begin{equation}\label{eq:Psi:s}
\Psi_{\pi_0}(K,E,\alpha)\le2\pi_0+s+6\pi_0\sqrt K\logp\left(\frac{1/\alpha-2\pi_0}{s}\right)+\frac{4\pi_0E}s.
\end{equation}
In particular, with $s^*=3\pi_0\sqrt K+\sqrt{9\pi_0^2K+4\pi_0E}$,
\begin{align*}
\Psi_{\pi_0}(K,E,\alpha)&\le2\pi_0+2s^*-6\pi_0\sqrt K+6\pi_0\sqrt K\,\logp\left(\frac{1/\alpha-2\pi_0}{s^*}\right)\\
&\le\pi_0\big(2+6\sqrt K\big)+4\sqrt{\pi_0E}+6\pi_0\sqrt K\,\logp\frac1{6\alpha\pi_0\sqrt K},
\end{align*}
where $\logp x=\max\{\log x,0\}$.
\end{enumerate}
\end{lemma}
\begin{proof}
\noindent\ref{lem:Psi:ii}. By \eqref{eq:def:Psi}, $\Psi_{\pi_0}(K,E,\alpha)=\int_0^{1/\alpha}\psi(t)\,dt$, where $\psi(t)=1$ for $t\le2\pi_0$ and $\psi(t)=\min\{1,g(t)\}$ for $t>2\pi_0$. For each $t$, $\psi(t)$ is nondecreasing in $K$ and in $E$. It is also nondecreasing in $\pi_0$, since $\pi_0/(t-2\pi_0)$ and $\pi_0/(t-2\pi_0)^2$ increase with $\pi_0$ for $t>2\pi_0$, and $\psi(t)=1$ once $\pi_0\ge t/2$. Finally $0\le\psi\le1$ and $\psi=1$ on $[0,2\pi_0]$, which gives the two bounds.

\noindent\ref{lem:Psi:iii}. For every $s>0$, bounding $\min\{1,g(t)\}$ by $1$ for $2\pi_0<t\le2\pi_0+s$ and by $g(t)$ for $t>2\pi_0+s$ in \eqref{eq:def:Psi} gives \eqref{eq:Psi:s}. Here $\logp$ accounts for the case $2\pi_0+s\ge1/\alpha$, in which the integral of $g$ is over an empty interval. Take $s=s^*$, the positive root of $s^2-6\pi_0\sqrt Ks-4\pi_0E=0$. Then $4\pi_0E/s^*=s^*-6\pi_0\sqrt K$, which gives the first inequality. For the second, $\sqrt{9\pi_0^2K+4\pi_0E}\le3\pi_0\sqrt K+2\sqrt{\pi_0E}$ gives $s^*\le6\pi_0\sqrt K+2\sqrt{\pi_0E}$, and $s^*\ge6\pi_0\sqrt K$ and $1/\alpha-2\pi_0\le1/\alpha$ give $\logp((1/\alpha-2\pi_0)/s^*)\le\logp(1/(6\alpha\pi_0\sqrt K))$. Proof is completed.
\end{proof}
For $\theta\in(0,1)$ define
\[
E_{\rho,m}(\alpha)=\frac{\Lambda_{\rho,m}(\alpha)}{\theta^2}+\frac{1}{(1-\theta)^2}\sum_{j=0}^{J-1}e_{\rho,m}(z_j),
\]
with $\Lambda_{\rho,m}(\alpha)$ and $e_{\rho,m}$ from Section~\ref{sec:prelim}. Recall that $b_k$ in $\Lambda_{\rho,m}(\alpha)$ is a supremum over $0\le z\le z_{J-1}$, where $z_{J-1}$ is the largest threshold.

\begin{proposition}\label{prop:master}
Let $\Sigma$ be an $m\times m$ correlation matrix and $\rho=\rho_\Sigma$. Let $\alpha\in(0,1/2]$, $m\ge2$ and $\theta\in(0,1)$. For one- and two-sided $p$-values and every $\mu\in\mathbb R^m$,
\begin{equation}\label{eq:initial:FDR}
\FDR(\alpha;m,\mu,\Sigma)\le\alpha\min\Big\{\pi_0H_m,\ \Psi_{\pi_0}\Big(\frac{K_\rho}{(1-\theta)^2},\,E_{\rho,m}(\alpha),\,\alpha\Big)\Big\},
\end{equation}
where $H_m=\sum_{k=1}^mk^{-1}\sim\log{m}$.
\end{proposition}
\begin{proof}
The first term is the bound of \citet{BY2001}; see (\ref{BH:original:bound}). For the second, we may assume $m_0\ge1$, since otherwise $\FDR=0$. Write $z_j=z_{w_j}$, $K'=K_\rho/(1-\theta)^2$ and $\bar H_k=m_0^{-1}\sum_{i\in\mathcal H_0}He_k(\xi_i)$.

\emph{Step 1: splitting the centered count.} By the definition of $Q_i$,
\[
\Ntil(w_j)-\pi_0a_j=m_0\sum_{k\in\mathcal S,\ k<k^*}\gamma_k(z_j)\bar H_k+\sum_{i\in\mathcal H_0}Q_i(w_j).
\]
Since $|\gamma_k(z_j)|/w_j=|c_k(z_j)|/(1-\Phi(z_j))\le b_k$, the first term is at most $m_0w_j\Lambda^*=\pi_0a_j\Lambda^*$ with $\Lambda^*=\sum_{k<k^*}b_k|\bar H_k|$, which does not depend on $j$. By Lemma~\ref{lem:hermite}\ref{lem:hermite:ii}, 
\[
\E\bar H_k^2=k!\,m_0^{-2}\sum_{i,i'\in\mathcal H_0}\Sigma_{ii'}^k\le k!(1/m_0+2G_{\rho,k}(m)/m_0),
\]
because $|\Sigma_{ii'}|\le\rho(|i-i'|)$ and each $i$ has at most two partners at each lag. Hence $\pi_0\E\bar H_k^2\le k!(1/m+2G_{\rho,k}(m)/m)$, and Minkowski's inequality gives $\pi_0\E(\Lambda^*)^2\le\Lambda_{\rho,m}(\alpha)$.

\emph{Step 2: a tail bound for $\Stil$.} Fix $t>2\pi_0$ and let $B_j=\{\Ntil(w_j)\ge\frac t2a_j\}$. On $B_j$, $\Ntil(w_j)-\pi_0a_j\ge\frac{t-2\pi_0}2a_j$. So $B\cup B_j^2$ holds, where
\[
B=\Big\{\pi_0\Lambda^*\ge\theta\frac{t-2\pi_0}2\Big\},\qquad B_j^2=\Big\{\sum_{i\in\mathcal H_0}Q_i(w_j)\ge(1-\theta)\frac{t-2\pi_0}2a_j\Big\}.
\]
Moreover $\Ntil(w_j)>0$ on $B_j$, so $B_j^1\triangleq\{\Ntil(w_j)\ge1\}$ holds. By Lemma~\ref{lem:dyadic},
\[
\{\Stil\ge t\}\subseteq\bigcup_{j=0}^{J-1}\big(B_j^1\cap(B\cup B_j^2)\big)\subseteq B\cup\bigcup_{j=0}^{J-1}\big(B_j^1\cap B_j^2\big),
\]
and therefore
\[
\P(\Stil\ge t)\le\P(B)+\sum_{j=0}^{J-1}\min\{\P(B_j^1),\P(B_j^2)\}.
\]
By Markov's inequality, $\P(B_j^1)\le\E\Ntil(w_j)=\pi_0a_j$. By Chebyshev's inequality and Step~1,
\[
\P(B)\le\frac{4\pi_0^2\E(\Lambda^*)^2}{\theta^2(t-2\pi_0)^2}\le\frac{4\pi_0\Lambda_{\rho,m}(\alpha)}{\theta^2(t-2\pi_0)^2}.
\]
By Chebyshev's inequality and Proposition~\ref{prop:var} with $I=\mathcal H_0$,
\[
\P(B_j^2)\le\frac{4(K_\rho m_0w_j+e_{\rho,m}(z_j)m_0mw_j^2)}{(1-\theta)^2(t-2\pi_0)^2a_j^2}=\frac b{a_j}+c_j,
\]
where $b=4\pi_0K'/(t-2\pi_0)^2$ and $c_j=4\pi_0e_{\rho,m}(z_j)/((1-\theta)^2(t-2\pi_0)^2)$. Since $\min\{x,y+c\}\le\min\{x,y\}+c$ for $c\ge0$,
\[
\min\{\P(B_j^1),\P(B_j^2)\}\le\min\Big\{\pi_0a_j,\frac b{a_j}\Big\}+c_j=\min\Big\{\pi_0a_j,\frac{\pi_0b}{\pi_0a_j}\Big\}+c_j.
\]
The numbers $\pi_0a_j$ halve from one level to the next, so Lemma~\ref{lem:geom} gives
\[
\sum_{j=0}^{J-1}\min\Big\{\pi_0a_j,\frac{\pi_0b}{\pi_0a_j}\Big\}\le3\sqrt{\pi_0b}=\frac{6\pi_0\sqrt{K'}}{t-2\pi_0}.
\]
By the definition of $E_{\rho,m}(\alpha)$,
\[
\P(B)+\sum_{j=0}^{J-1}c_j\le\frac{4\pi_0}{(t-2\pi_0)^2}\Big(\frac{\Lambda_{\rho,m}(\alpha)}{\theta^2}+\frac{1}{(1-\theta)^2}\sum_{j=0}^{J-1}e_{\rho,m}(z_j)\Big)=\frac{4\pi_0E_{\rho,m}(\alpha)}{(t-2\pi_0)^2}.
\]
Therefore, for every $t>2\pi_0$,
\[
\P(\Stil\ge t)\le\frac{6\pi_0\sqrt{K'}}{t-2\pi_0}+\frac{4\pi_0E_{\rho,m}(\alpha)}{(t-2\pi_0)^2}.
\]

\emph{Step 3: integration.} Since $2\pi_0\le1/\alpha$, we split the integral in Lemma~\ref{lem:reduction} at $2\pi_0$. We use $\P(\Stil\ge t)\le1$ for $t\le2\pi_0$, and both this bound and Step~2 for $t>2\pi_0$. By \eqref{eq:def:Psi},
\begin{align*}
\frac{\FDR(\alpha;m,\mu,\Sigma)}{\alpha}&\le\int_0^{2\pi_0}1\,dt+\int_{2\pi_0}^{1/\alpha}\P(\Stil\ge t)\,dt\\
&\le2\pi_0+\int_{2\pi_0}^{1/\alpha}\min\Big\{1,\frac{6\pi_0\sqrt{K'}}{t-2\pi_0}+\frac{4\pi_0E_{\rho,m}(\alpha)}{(t-2\pi_0)^2}\Big\}dt=\Psi_{\pi_0}\big(K',E_{\rho,m}(\alpha),\alpha\big).
\end{align*}
Proof is completed.
\end{proof}

\subsection{Proofs of Theorem~\ref{thm:main}}\label{sec:proofmain}
\emph{Monotonicity in the correlations.} The quantities $G_{\rho,k}(m)$, $e_{\rho,m}(z)$, $\Lambda_{\rho,m}(\alpha)$ and $E_{\rho,m}(\alpha)$ depend on $\rho$ only through $|\rho(l)|$, $l\ge1$, and are nondecreasing in each $|\rho(l)|$. The same holds for $K_\rho$. Indeed, $R(\cdot,z)$ is increasing by Lemma~\ref{lem:hermite}\ref{lem:hermite:iii}, so $\tau(r,z)$ is nondecreasing in $|r|$, and increasing $|\rho(l)|$ can only add lags to the sum in $K_\rho$. These quantities are well defined at the envelope $\rho_{\epsilon,C}$, although the matrix with entries $\rho_{\epsilon,C}(|i-j|)$ need not be positive semidefinite.
 
\emph{Proof of \ref{thm:main:i}.} For $m=1$, we have $H_1=1$ and $F_{\epsilon,C,1}(\alpha;\pi_0)=\pi_0$, so \eqref{our:FDR:bound} reads $\FDR\le\alpha\pi_0$, which is the bound \eqref{BH:original:bound}. Let $m\ge2$, $\theta\in(0,1)$, $\Sigma\in\mathcal C_m(\epsilon,C)$ and $\rho=\rho_\Sigma$. Since $\rho\le\rho_{\epsilon,C}$ pointwise, the monotonicity above gives $K_\rho\le K_{\epsilon,C}$, $\Lambda_{\rho,m}(\alpha)\le\Lambda_{\epsilon,C,m}(\alpha)$ and $\sum_je_{\rho,m}(z_j)\le E_{\epsilon,C,m}(\alpha)$. Hence 
\[
E_{\rho,m}(\alpha)\le E_m(\theta):=\Lambda_{\epsilon,C,m}(\alpha)/\theta^2+E_{\epsilon,C,m}(\alpha)/(1-\theta)^2.
\]
As $\Psi_{\pi_0}$ is nondecreasing in its first two arguments (Lemma~\ref{lem:Psi}\ref{lem:Psi:ii}), Proposition~\ref{prop:master} yields, for every $\mu\in\bbR^m$,
\begin{equation}\label{FDR:m:Psi}
\FDR(\alpha;m,\mu,\Sigma)\le\alpha\min\Big\{\pi_0H_m,\ \Psi_{\pi_0}\Big(\frac{K_{\epsilon,C}}{(1-\theta)^2},E_m(\theta),\alpha\Big)\Big\}.
\end{equation}
Taking the infimum over $\theta\in(0,1)$ gives \eqref{our:FDR:bound}. If $k^*=\min\mathcal S$, then $\Lambda_{\epsilon,C,m}(\alpha)=0$. The arguments of $\Psi_{\pi_0}$ in \eqref{eq:Fm} decrease to $(K_{\epsilon,C},E_{\epsilon,C,m}(\alpha))$ as $\theta\downarrow0$. Since $\Psi_{\pi_0}$ is nondecreasing and continuous in $(K,E)$, the infimum equals $\Psi_{\pi_0}(K_{\epsilon,C},E_{\epsilon,C,m}(\alpha),\alpha)$.
 
\emph{Proof of \eqref{eq:supF}.} If $\pi_0=0$, then $\Psi_0=0$, so $F_{\epsilon,C,m}(\alpha;0)=0$ and \eqref{eq:supF} holds. So we only consider $\pi_0\in(0,1]$.
If $k^*=\min\mathcal S$, then 
$\Lambda_{\epsilon,C,m}(\alpha)=0$, leading to
$F_{\epsilon,C,m}(\alpha;\pi_0)=\Psi_{\pi_0}(K_{\epsilon,C},E_{\epsilon,C,m}(\alpha),\alpha)$. Otherwise, taking $\theta=1/2$ in \eqref{eq:Fm} gives
\[
F_{\epsilon,C,m}(\alpha;\pi_0)\le\Psi_{\pi_0}\big(4K_{\epsilon,C},4\Lambda_{\epsilon,C,m}(\alpha)+4E_{\epsilon,C,m}(\alpha),\alpha\big).
\]
By Lemma~\ref{lem:Psi}\ref{lem:Psi:ii}, it therefore suffices to bound $K_{\epsilon,C}$, $\Lambda_{\epsilon,C,m}(\alpha)$ and $E_{\epsilon,C,m}(\alpha)$ uniformly in $m\ge2$. We do this in Steps~1--3 and conclude in Step~4.

\emph{Reduction to $\epsilon=1/k^*$.} Since $\log(l+e)\ge1$, the envelope $\rho_{\epsilon,C}$ is nonincreasing in $\epsilon$. As $\epsilon\ge1/k^*$, we have $\rho_{\epsilon,C}\le\rho_{1/k^*,C}$. Let $\rho=\rho_{1/k^*,C}$ and $p=1+1/k^*$, so that $|\rho(l)|\le C\log^{-p}(l+e)$ for $l\ge1$ and $k^*p=k^*+1$. By the monotonicity above,
\[
K_{\epsilon,C}\le K_\rho,\qquad \Lambda_{\epsilon,C,m}(\alpha)\le\Lambda_{\rho,m}(\alpha),\qquad E_{\epsilon,C,m}(\alpha)\le\sum_{j=0}^{J-1}e_{\rho,m}(z_j).
\]
We bound the three right-hand sides. Recall $J=\lceil\log_2m\rceil$, and let $m\ge2$.

\emph{Step 1: bound on $K_{\epsilon,C}$.} For $z\le\sqrt2$ we have $r_z=1/2$, and each lag with $|\rho(l)|>1/2$ contributes $2$ to the sum in $K_\rho$. For $z>\sqrt2$ we have $r_z=z^{-2}<1/2$. Then each lag with $|\rho(l)|>1/2$ still contributes $2$, and each lag with $z^{-2}<|\rho(l)|\le1/2$ contributes at most $2(1+z)e^{-z^2/6}/\sqrt{6\pi}$ by Lemma~\ref{lem:mills}\ref{lem:mills:ii}. Since $|\rho(l)|>c$ implies $l<\exp\{(C/c)^{1/p}\}-e$, there are at most $(e^{(2C)^{1/p}}-e)_+$ lags of the first kind and at most $(\exp\{(Cz^2)^{1/p}\}-e)_+$ lags of the second kind. Hence
\begin{equation}
K_{\epsilon,C}\le\bar K_{C,k^*}:=1+2\big(e^{(2C)^{1/p}}-e\big)_++\frac{2}{\sqrt{6\pi}}\sup_{z>\sqrt2}(1+z)e^{-z^2/6}\Big(\exp\big\{(Cz^2)^{1/p}\big\}-e\Big)_+.
\end{equation}
Since $p>1$, we have $(Cz^2)^{1/p}=o(z^2)$ as $z\to\infty$, so $\bar K_{C,k^*}<\infty$. It depends only on $C$ and $k^*$, through $p=1+1/k^*$.

\emph{Two elementary facts.} Since $m\ge2$, we have $1\le J\le(2/\log2)\log m$. For $q>0$ let $g_q(m)=m^{-1}\sum_{l=1}^{m-1}\log^{-q}(l+e)$, so that $G_{\rho,k}(m)/m\le C^kg_{kp}(m)$. At most $\sqrt m$ lags satisfy $l<\sqrt m$, each with a term at most $1$, and $\log(l+e)\ge\frac12\log m$ for $l\ge\sqrt m$. Hence $g_q(m)\le m^{-1/2}+2^q\log^{-q}m$. Using $\sup_{x>0}x^ae^{-x/c}=(ac/e)^a$ with $x=\log m$, we obtain
\begin{align}
J^{k^*+1}g_{k^*+1}(m)&\le M_{k^*}:=\Big(\frac2{\log2}\Big)^{k^*+1}\Big[\Big(\frac{2k^*+2}e\Big)^{k^*+1}+2^{k^*+1}\Big],\label{eq:M}\\
J^k\Big(\frac1m+2C^kg_{kp}(m)\Big)&\le N_k:=\Big(\frac2{\log2}\Big)^k\Big[\Big(\frac ke\Big)^k+2C^k\Big(\Big(\frac{2k}e\Big)^k+\frac{2^{kp}}{\log2}\Big)\Big],\quad k\in\mathcal S,\ k<k^*.\label{eq:N}
\end{align}
In \eqref{eq:N} we used $\log^{k-kp}m=\log^{-k/k^*}m\le1/\log2$. Moreover, by Lemma~\ref{lem:mills}\ref{lem:mills:i}, $z_w^2\le2\log(1/w)$ in both cases, so
\begin{equation}\label{eq:zj}
z_j^2\le2\log(1/\alpha)+2j\log2,\qquad 0\le j\le J-1.
\end{equation}

\emph{Step 2: bound on $E_{\epsilon,C,m}(\alpha)$.} By Lemma~\ref{lem:mills}\ref{lem:mills:i}, $(\phi(z)/(1-\Phi(z)))^2\le z^2+2$. By Lemma~\ref{lem:A}\ref{lem:A:iii}, $\mathsf A(z)\le a_1(z^2+2)^{k^*-1}$. In the two-sided case with $k^*=2$, Lemma~\ref{lem:A}\ref{lem:A:ii} gives $\mathsf A(z)\le a_2(z^2+2)$. Let $a_*=a_2$ in this case and $a_*=a_1$ otherwise. Then $e_{\rho,m}(z)\le2a_*C^{k^*}(z^2+2)^{k^*}g_{k^*+1}(m)$. By \eqref{eq:zj}, $z_j^2+2\le2(1+\log(1/\alpha))(1+j)$. Since $\sum_{j=0}^{J-1}(1+j)^{k^*}\le J^{k^*+1}$, \eqref{eq:M} gives
\begin{equation}
E_{\epsilon,C,m}(\alpha)\le\bar E_{C,k^*}\big(1+\log(1/\alpha)\big)^{k^*},\qquad \bar E_{C,k^*}:=2^{k^*+1}a_*C^{k^*}M_{k^*}.
\end{equation}

\emph{Step 3: bound on $\Lambda_{\epsilon,C,m}(\alpha)$.} If $k^*=\min\mathcal S$, then $\Lambda_{\epsilon,C,m}(\alpha)=0$. Otherwise let $k\in\mathcal S$ with $k<k^*$, and write $Z=z_{J-1}$. By Lemma~\ref{lem:mills}\ref{lem:mills:i},\ref{lem:mills:iii}, $\phi(z)/(1-\Phi(z))\le h(z)$ and $|He_{k-1}(z)|\le(z+\sqrt k)^{k-1}$, and both are nondecreasing in $z$. Using $h(Z)^2\le Z^2+2$ and $(Z+\sqrt k)^2\le2(Z^2+k)$, we get $b_k^2\,k!\le2^{k-1}(Z^2+k+1)^k/k!$. By \eqref{eq:zj} and $J\ge1$, $Z^2+k+1\le(k+1)(1+\log(1/\alpha))+2J\le2(k+1)(1+\log(1/\alpha))J$. Together with \eqref{eq:N},
\[
b_k^2\,k!\Big(\frac1m+\frac{2G_{\rho,k}(m)}m\Big)\le\frac{2^{k-1}(2k+2)^k}{k!}N_k\big(1+\log(1/\alpha)\big)^{k^*}.
\]
Let $N=\#\{k\in\mathcal S:k<k^*\}$. By the Cauchy--Schwarz inequality,
\begin{equation}
\Lambda_{\epsilon,C,m}(\alpha)\le\bar\Lambda_{C,k^*}\big(1+\log(1/\alpha)\big)^{k^*},\qquad \bar\Lambda_{C,k^*}:=N\sum_{k\in\mathcal S,\ k<k^*}\frac{2^{k-1}(2k+2)^k}{k!}N_k.
\end{equation}
If $k^*=\min\mathcal S$, then $N=0$ and $\Lambda_{\epsilon,C,m}(\alpha)=0$, so \eqref{eq:Lbar} holds in both cases.

\emph{Step 4: summary.} By Steps~1--3 and Lemma~\ref{lem:Psi}\ref{lem:Psi:ii}, in both cases $k^*=\min\mathcal S$ and $k^*>\min\mathcal S$, for every $m\ge2$,
\[
F_{\epsilon,C,m}(\alpha;\pi_0)\le\Psi_{\pi_0}\Big(4\bar K_{C,k^*},\ 4\big(\bar\Lambda_{C,k^*}+\bar E_{C,k^*}\big)\big(1+\log(1/\alpha)\big)^{k^*},\ \alpha\Big),
\]
with $\bar K_{C,k^*}$, $\bar\Lambda_{C,k^*}$ and $\bar E_{C,k^*}$ given in \eqref{eq:Kbar}, \eqref{eq:Lbar} and \eqref{eq:Ebar}. Apply the second inequality of Lemma~\ref{lem:Psi}\ref{lem:Psi:iii} with $K=4\bar K_{C,k^*}$ and $E=4(\bar\Lambda_{C,k^*}+\bar E_{C,k^*})(1+\log(1/\alpha))^{k^*}$. Then $6\pi_0\sqrt K=12\pi_0\sqrt{\bar K_{C,k^*}}$ and $4\sqrt{\pi_0E}=8\sqrt{\pi_0(\bar\Lambda_{C,k^*}+\bar E_{C,k^*})}\,(1+\log(1/\alpha))^{k^*/2}$. Since $\logp(xy)\le\logp x+\logp y$, we have $\logp(1/(12\alpha\pi_0\sqrt{\bar K_{C,k^*}}))\le\log(1/\alpha)+\logp(1/(12\pi_0\sqrt{\bar K_{C,k^*}}))$. Hence
\[
\begin{aligned}
\sup_{m\ge2}F_{\epsilon,C,m}(\alpha;\pi_0)\le{}&\pi_0\big(2+12\sqrt{\bar K_{C,k^*}}\big)+12\pi_0\sqrt{\bar K_{C,k^*}}\Big(\log(1/\alpha)+\logp\frac1{12\pi_0\sqrt{\bar K_{C,k^*}}}\Big)\\
&+8\sqrt{\pi_0\big(\bar\Lambda_{C,k^*}+\bar E_{C,k^*}\big)}\,\big(1+\log(1/\alpha)\big)^{k^*/2}.
\end{aligned}
\]
Finally, let $M=\max\{1,k^*/2\}$ and $x=\log(1/\alpha)\ge0$. Since $1+x\ge1$ and $k^*/2\le M$, convexity gives $(1+x)^{k^*/2}\le(1+x)^M\le2^{M-1}(1+x^M)$. Together with $x\le1+x^M$, this gives \eqref{eq:supF}. This completes the proof of Theorem~\ref{thm:main}. \qed

\bibliographystyle{plainnat}
\bibliography{refs}

\end{document}